\documentclass[11pt]{article}

\usepackage{import,mathrsfs,float,pgfplots,mhchem,
bbm,graphicx,subcaption,listings,xcolor,tcolorbox}\usepackage{etoolbox,booktabs}

\usepackage{caption}
\usepackage[acronym,toc,sort=use,nonumberlist]{glossaries}

\preto\printglossary{\glsaddall}
\makeglossaries

\lstdefinestyle{mathematica}{
    language=Mathematica,
    basicstyle=\small\ttfamily,
    backgroundcolor=\color{gray!10},
    frame=single,
    breaklines=true,
    commentstyle=\color{green!60!black},
    keywordstyle=\color{blue},
    stringstyle=\color{red},
    showstringspaces=false,
    tabsize=2
}
\usepackage[final,backref]{hyperref}
\usepackage{tikz-cd}
\usepackage{lineno}
\tikzcdset{every label/.append style = {font = \small}}
\usepackage{tikz}
\usetikzlibrary{graphs,decorations.pathmorphing,decorations.markings,arrows.meta,
automata,positioning}
\usepackage{subcaption}
\usepackage{nomencl}
\makenomenclature

\usepackage{makecell}
\usepackage{amssymb, amsfonts,mathtools,bbm,eurosym,graphicx,latexsym, amsmath,amsthm,times,nccmath,url}
\usepackage{listings,verbatim,cases,hyperref,tabularx}
\usepackage{graphicx,float,fancybox}
\usepackage{epsfig,epstopdf,tikz}
\usepackage[sans]{dsfont}
\usetikzlibrary{positioning,arrows,arrows.meta,calc}
\lstdefinestyle{mathematica}{
 language=Mathematica,
 basicstyle=\small\ttfamily,
 backgroundcolor=\color{gray!10},
 frame=single,
 breaklines=true,
 commentstyle=\color{green!60!black},
 keywordstyle=\color{blue},
 stringstyle=\color{red},
 showstringspaces=false,
 tabsize=2
}
\newtcolorbox{keyresult}[1][]{
 colback=blue!5!white,
 colframe=blue!75!black,
 title=#1,
 fonttitle=\bfseries
}

\newtcolorbox{examplebox}[1][]{
 colback=green!5!white,
 colframe=green!75!black,
 title=#1,
 fonttitle=\bfseries
}

\newtheorem{theorem}{Theorem}
\newtheorem{definition}{Definition}
\newtheorem{open}{Problem}
\newtheorem{example}{Example}

\newcounter{tabcounter}
\newcounter{defcounter}
\newcounter{opcounter}
\newcounter{thmcounter}
\newcounter{excounter}
\newcounter{ascounter}
\newcounter{cjcounter}

\newglossarystyle{unifiedstyle}{%
  \setglossarystyle{list}%
  \renewcommand{\glossentry}[2]{%
    \edef\glscursec{\glsentryuseriii{##1}}%
    \ifx\glscursec\glslastsec\else
      \item[]\textbf{\S\,\glscursec}%
      \global\let\glslastsec\glscursec
    \fi
    \item[\glsentryuserii{##1} \glsentryuseri{##1}:] ``\glossentryname{##1}'' \glsentrydesc{##1}%
  }%
}

\makeatletter
\newcommand{\beXa}[1][]{%
  \def\@temparg{#1}%
  \ifx\@temparg\@empty
    \begin{example}%
  \else
    \refstepcounter{excounter}
    \begin{example}[#1]%
    \label{e:\theexcounter}%
    \expanded{%
      \noexpand\newglossaryentry{e\theexcounter}{%
        name={\unexpanded{#1}},%
        sort={\theexcounter},%
        description={p.~\noexpand\pageref{e:\theexcounter}},%
        user1={\theexcounter},%
        user2={Example},%
        user3={\thesection}%
      }%
    }%
  \fi
}
\newcommand{\eeXa}{\end{example}}

\newcommand{\beD}[1][]{%
  \def\@temparg{#1}%
  \ifx\@temparg\@empty
    \begin{definition}%
  \else
    \refstepcounter{defcounter}%
    \begin{definition}[#1]%
    \label{d:\thedefcounter}%
    \expanded{%
      \noexpand\newglossaryentry{d\thedefcounter}{%
        name={\unexpanded{#1}},%
        sort={\thedefcounter},%
        description={p.~\noexpand\pageref{d:\thedefcounter}},%
        user1={\thedefcounter},%
        user2={Definition},%
        user3={\thesection}%
      }%
    }%
  \fi
}
\newcommand{\eeD}{\end{definition}}

\newcommand{\beT}[1][]{%
  \def\@temparg{#1}%
  \ifx\@temparg\@empty
    \begin{theorem}%
  \else
    \refstepcounter{thmcounter}%
    \begin{theorem}[#1]%
    \label{t:\thethmcounter}%
    \expanded{%
      \noexpand\newglossaryentry{t\thethmcounter}{%
        name={\unexpanded{#1}},%
        sort={\thethmcounter},%
        description={p.~\noexpand\pageref{t:\thethmcounter}},%
        user1={\thethmcounter},%
        user2={Theorem},%
        user3={\thesection}%
      }%
    }%
  \fi
}
\newcommand{\eeT}{\end{theorem}}
\newcommand{\beO}[1][]{%
  \def\@temparg{#1}%
  \ifx\@temparg\@empty
    \begin{open}%
  \else
    \refstepcounter{opcounter}%
    \begin{open}[#1]%
    \label{op:\theopcounter}%
    \expanded{%
      \noexpand\newglossaryentry{o\theopcounter}{%
        name={\unexpanded{#1}},%
        sort={\theopcounter},%
        description={p.~\noexpand\pageref{op:\theopcounter}},%
        user1={\theopcounter},%
        user2={Problem},%
        user3={\thesection}%
      }%
    }%
  \fi
}
\newcommand{\eeO}{\end{open}}
\newcommand{\beAs}[1][]{%
  \def\@temparg{#1}%
  \ifx\@temparg\@empty
    \begin{assumption}%
  \else
    \refstepcounter{ascounter}%
    \begin{assumption}[#1]%
    \label{as:\theascounter}%
    \expanded{%
      \noexpand\newglossaryentry{as\theascounter}{%
        name={\unexpanded{#1}},%
        sort={\theascounter},%
        description={p.~\noexpand\pageref{as:\theascounter}},%
        user1={\theascounter},%
        user2={Assumption},%
        user3={\thesection}%
      }%
    }%
  \fi
}
\newcommand{\eeAs}{\end{assumption}}
\newcommand{\beCj}[1][]{%
  \def\@temparg{#1}%
  \ifx\@temparg\@empty
    \begin{conjecture}%
  \else
    \refstepcounter{cjcounter}%
    \begin{conjecture}[#1]%
    \label{cj:\thecjcounter}%
    \expanded{%
      \noexpand\newglossaryentry{cj\thecjcounter}{%
        name={\unexpanded{#1}},%
        sort={\thecjcounter},%
        description={p.~\noexpand\pageref{cj:\thecjcounter}},%
        user1={\thecjcounter},%
        user2={Conjecture},%
        user3={\thesection}%
      }%
    }%
  \fi
}
\newcommand{\eeCj}{\end{conjecture}}

\newcommand{\betA}[1][]{%
  \def\@temparg{#1}%
  \ifx\@temparg\@empty
    \begin{table}[htbp]%
  \else
    \refstepcounter{tabcounter}%
    \begin{table}[htbp]%

    \label{A:\thetabcounter}%

    \expanded{%
      \noexpand\newglossaryentry{A\thetabcounter}{%
        name={\unexpanded{#1}},%
        sort={\thetabcounter},%
        description={p.~\noexpand\pageref{A:\thetabcounter}},%
        user1={\thetabcounter},%
        user2={Table},%
        user3={\thesection}%
      }%
    }%
  \fi
}

\newcommand{\eetA}{\end{table}}
\makeatother

\newtheorem{lemma}{Lemma}
\newtheorem{proposition}{Proposition}
\newtheorem{corollary}{Corollary}
\newtheorem{conjecture}{Conjecture}
\newtheorem{remark}{Remark}
\newtheorem{question}{Question}
\newtheorem{assumption}{Assumption}
\newtheorem{hypothesis}{Hypothesis}

\def\beL{\begin{lemma}}\def\eeL{\end{lemma}}
\def\beP{\begin{proposition}}\def\eeP{\end{proposition}}
\def\beC{\begin{corollary}}\def\eeC{\end{corollary}}
\def\beR{\begin{remark}}\def\eeR{\end{remark}}
\def\beQ{\begin{question}}\def\eeQ{\end{question}}
\def\beA{\begin{assumption}}\def\eeA{\end{assumption}}

\definecolor{funccolor}{RGB}{25,25,112}
\definecolor{desccolor}{RGB}{64,64,64}

\def\bep{\begin{pmatrix}}\def\eep{\end{pmatrix}}
\def\bev{\begin{vmatrix}}\def\eev{\end{vmatrix}}
\def\bea{\begin{eqnarray*}}\def\eea{\end{eqnarray*}}
\def\bc{\begin{cases}}\def\ec{\end{cases}}
\def\BEN{\begin{enumerate}}\def\EEN{\end{enumerate}}
\def\BI{\begin{itemize}}\def\EI{\end{itemize}}

\newcommand{\be}[1]{\begin{equation}\label{#1}}
\newcommand{\ee}{\end{equation}}
\newcommand{\beq}{\begin{eqnarray}}
\def\eeq{\end{eqnarray}}
\def\eqr{\eqref}\def\fr{\frac}\def\lbl{\label}
\def\Lra{\Longrightarrow}\def\cNGM{\cite{Diek,Van,Van08}}
\def\Eq{\Leftrightarrow}\def\b{\beta}

\newcommand{\R}{\mathbb{R}}  

\newcommand{\Z}{\mathbb{Z}}

\def\Se{\bar{S}} \def\Ie{\bar{I}}

\def\ba{\bff \alpha}

\def\eps{\varepsilon}
\def\f{\varphi}  
\def\Ga{\Gamma}\def\G{\Gamma}
\def\la{\lambda}\def\La{\Lambda}

\def\Gap{\Gamma^{+}} \def\Gam{\Gamma^{-}}

\def\T{\widetilde}

\newcommand{\bff}[1]{{\mbox{\boldmath$#1$}}}

\def\x{\boldsymbol{x}}\def\xe{\bar{\x}}
\newcommand\y{\boldsymbol{y}}\def\ye{\bar{\y}}

\def\vbe{\vec {\bff \beta}}

\def\v1{\vec {\bff 1}}

\def\mR{{\mathcal R}}

\DeclareMathOperator{\diag}{Diag}

\newcommand\op{\operatorname}
\newcommand{\supp}{\op{supp}}

\providecommand{\pr}[1]{\left(#1\right)}
\providecommand{\pp}[1]{\left[#1\right]}

\newcommand\blue[1]{\textcolor{blue}{#1}}
\newcommand{\red}{\textcolor[rgb]{1.00,0.00,0.00}}

\long\def\symbolfootnote[#1]#2{%
\begingroup
\def\thefootnote{\fnsymbol{footnote}}\footnote[#1]{#2}%
\endgroup}

\def\and{antisymmetric}

\def\evec{eigenvector}
\def\elt{element}

\def\ie{i.e. }
\def\im{\item}

\def\mbw{may be written as}

\def\nne{non-negative}

\def\rep{representation}
\def\resp{respectively}
\def\saty{satisfy}
\def\ssec{\subsection}

\def\wrt{with respect to}

\newcommand\CRN{chemical reaction network}

\def\brn{basic reproduction number}
\def\com{compartment}
\def\DFE{disease free equilibrium}

\def\EE{endemic point}

\def\NGM{next generation matrix}

\def\jin{Jacobian matrix of the invasion vector field with respect to the invasion variables}

\def\bfp{boundary fixed point}

\def\GAS{globally asymptotically stable}

\def\MH{Metzler-Hurwitz}\def\RH{Routh-Hurwitz}

\def\Lf{Lyapunov function}

\def\ABa{$(A,\bbe,\ba,\varphi)$ balanced bilinear model}
\def\ABP{
\\ balanced bilinear model}

\def\gLV{generalized Lotka-Volterra model}

\def\regS{{\bf regular splitting}}

 \def\ABP{$(A,\bbe,P,\f)$ balanced bilinear model}
\def\ABa{
rank-one balanced bilinear model}
\def\bbe{B
} \def\vbe{\vec {\bff \beta}
}
 
\newcommand{\tws}{\tilde{w}^{\sigma}}
 
\def\bw{{\mathbf w}}

\def\y{S} \def\x{I} \def\ye{\bar{S}}\def\xe{\bar{I}}
\newtcolorbox{proofbox}[1][]{
    colback=yellow!5!white,
    colframe=orange!75!black,
    title=#1,
    fonttitle=\bfseries
}

\title{A Perron-Frobenius  strong  threshold theorem for \ABP s,
and the role of left and right Perron \evec s in mathematical epidemiology\footnote{This work was supported in part by the National Science and Technology Council of Taiwan [grant numbers: 113-2222-E-110-002-MY3 and 114-2218-E-007-011-].}}

\author{Florin Avram$^{1}$,
Rim Adenane$^{2}$, Andrei-Dan Halanay$^{3}$, Andras Horvath$^{4}$ and Sei Zhen Khong$^{5}$
}

\begin{document}

\maketitle
\begin{center}	

$^{1}$ Laboratoire de Math\'{e}matiques Appliqu\'{e}es, Universit\'{e} de Pau, 64000, Pau, France; avramf3@gmail.com \\

$^{2}$ Laboratoire d'Analyse, G\'{e}om\'{e}trie et Applications, d\'epartement des Math\'ematiques, Universit\'e Ibn-Tofail, 14000,  Kenitra, Maroc; rim.adenane@uit.ac.ma \\

$^{3}$ Department of Mathematics and Computer Science, University of
Bucharest, Bucharest, RO-010014, Romania; halanay@fmi.unibuc.ro\\

 $^{4}$ Department of Computer Science, University of Turin,  10124, Italy; horvath@di.unito.it \\

$^{5}$ Department of Electrical Engineering, National Sun Yat-sen University, Kaohsiung, 80424, Taiwan; szkhong@mail.nsysu.edu.tw.

\end{center}

\begin{abstract}

This paper started as a review of  seven results pertaining to a family of bilinear models with rank one NGM introduced by Fall, Iggidr, Sallet and Bonzi, which utilize the explicit eigenvectors of the NGM to compute the unique endemic equilibrium (EE), and Lyapunov functions at both the disease free equilibrium (DFE) and EE, and of results of Shuai and Van den Driessche (2013), which essentially deal with the same ``DFE -EE stability exchange" in the non-rank one case, when the eigenvectors  are not explicit. Recently, these results were complemented  by Earn and McCluskey (2025), who proved as well a ``strong threshold theorem", namely that when the DFE is unstable, a second equilibrium which is globally asymptotically stable must exist.

Below, we obtain  in Theorem \ref{thm:TK_DFE} some results that extend beyond rank one. For example, a nontrivial positive equilibrium exists if and only if the spectral equation
\(\rho(\widetilde K(S))=1\),  admits a strictly positive solution.

Also, we showed that Bonzi-Iggidr-Sallet bilinear models with rank one NGM
may be classified in two classes, with slight variations in the eigenvector formulas, and that extensions in the presence of feedback from infectious to susceptible are possible.

Another takeout from the previous works, which we clarify  in a revisit of the seven results in the rank one case, is that Lyapunov functions for both the DFE and the EE may be constructed using as weights the left Perron eigenvector $\pi(S)$ of the \NGM\ (NGM) $K=F(S) V^{-1}$, where $S$ denote all the non-infectious variables, and when a positive ODE leaves  a siphon face, it does so along the right Perron eigenvector $w(S)$ of $\T K= V^{-1}F(S)$. The question of whether this continues to be true beyond rank one is explored in  ongoing work.
\end{abstract}

\textbf{Keywords:}  positive ODE;  disease free equilibrium; endemic equilibrium; Metzler matrices; next-generation matrix; regular splitting; Perron–Frobenius eigenvectors;  balanced bilinear models; WAIFW matrix;  Kirchoff's matrix-tree theorem.

\section{Introduction}
This paper started as a review of several recent results in
mathematical epidemiology (ME).  One  is an elegant proof of the celebrated next generation matrix  (NGM) method, which was shown recently \cite{AAH26,AH26}
to be an outcome  of  the fact that the siphon property used in chemical reaction networks (CRN) forces the Jacobian at any ``reasonable" boundary equilibrium to have a triangular block form. If in addition   the transversal-block is a Metzler matrix, then stability is amenable to next-generation analysis (which is analog to the \regS\ theory \cite{Varga} in the numerical analysis literature).

Another  result, due to Earn and McCluskey \cite{EarnMc}, is  that a general class of bilinear ME models,  first featured in papers like \cite{IggidrCEP,Fall,Bonzi}, have   an explicit \Lf\ which guarantees that the unique strictly positive ``endemic equilibrium" (EE)  is \GAS\ (GAS).

Note also a determinant law \cite{AAH22,AABGH}, which states that the sum of the determinants of the Jacobians at the DFE and EE equals 0, and had been proved for a SAIRS model only (but conjectured to hold more widely, together with    other explicit formulas/laws concerning both the EE and the DFE,  which hold for certain rank one models, and go back to \cite{IggidrCEP,Fall,Bonzi}).

While writing the review, we completed it with various remarks not stated explicitly in our sources, like the fact  that there are two classes of \ABP s,     and while five of the ``laws" hold for both of them, sometimes with small variations,  the determinant law and the \Lf\ at the endemic point hold only when there is just one susceptible  class, in which the two classes of \ABP s coincide.


\subsection{Positive dynamical systems}
A dynamical system is said to be positive
if it leaves invariant the positive orthant \cite{rantzer2015scalable}.

The biological interaction networks (BIN) sciences, which include
 population dynamics,
  ecology,
 mathematical epidemiology (ME),
  viro-immunology,
  \CRN (CRN), etc,
are all concerned with positive dynamical systems, and
 share common objects of study, such as multi-stationarity, local and global stability, bifurcations,  persistence, permanence, extinction, global stability and the construction of Lyapunov functions.   For interesting historical references on BIN, as well as recent results on persistence and global stability, the reader is referred to the CRN papers~\cite{CracNazPan,GopMilShiu,LyapPDE,Polly,craciun2019polynomial,
 CracDes,XuGAS}.

 While studying  quite similar problems, the BIN fields put emphasis on distinct  aspects.

 The ME literature focuses  often on non-polynomial models, on the effect of delays, and that of replacing ODEs by IDEs and PIDEs, which fit more realistically epidemiological data.

 The CRN literature focuses a lot on ``structural results", which are independent of the ODE rates, as long as these are smooth and monotone; it uses often graphs whose edges are labeled by the rates,  where the exact form of the rates may be unknown.

 There are very few results which are well-known and used equally in both fields (one exception being the works of Li and  Muldowney \cite{LiMul,LiMulGeo,LiMulSec}), and our hope is to bring the two fields closer to each other. 

Below,  we  review  and advance a chapter of the  ME literature, which could be called the ``laws of the DFE to EE stability relay", but also use on purpose occasionally the language  of  CRNT. While the theoretical results here do not rely on CRNT,
we also offer two example notebooks,  using our  Mathematica package Epid-CRN,
available at \url{https://github.com/florinav/EpidCRNmodels}, and readers willing to take advantage
of a ``hands on approach" will need some rudiments of CRNT.

\subsection{The (RHS,var), and the (RN,rts,var) mass-action representation of    positive ODEs}
{For parameterizing ME or CRN ODE models},
   there are three natural choices:

 \BEN \im
 The parametrization   used traditionally for  all ODE models   is $X'=RHS(X)$. A model is thus  a pair (RHS, var=X), where RHS denotes the set of formulas for the derivatives, and var are the  variables whose rates of change we study.

  \im Chemical reaction networks theory starts by breaking the RHS as \be{crnsys}X'=RHS(X)=\G .{\bf R}(X).\ee

    A model is now defined as a triple $(\G, {\bf R}(X), X)$.  The RHS has thus been broken into two parts:
  \BEN  \im
   $\G$ is the ``stoichiometric matrix'' (SM), whose columns  represent directions in which several species/\com s change simultaneously;
   it is also called the ``structure" of the model.
   \im  ${\bf R}(X)$ is the  vector of rates of change associated to each direction (assumed all to be \nne), also known as kinetics. They are separated from  $\G$, since they are less certain; this fact, well accepted in CRNT,  is equally true in ME.
    \EEN
\beXa[stoichiometric matrix of SIRS]
For example, the  simplest
 SIRS ODE without inflows and outflows (=demography) is defined by the  triple $(\G, {\bf R}, X)$:
  \be{SIRc}
X'=\bep
 s' \\  i'    \\  r'\eep
 =  \bep - 1&0&1&-1 \\
  1&- 1&0&0\\
  0&1&- 1&1
\eep \bep
\b s  i    \\ \gamma_i i    \\ \gamma_r r \\ \gamma_s s  \eep:=\Gamma {\bf R}(X).
\ee

\eeXa
Note this representation is also a natural first step towards defining an associated CTMC model (continuous-time Markov chain) on the integers.

We formalize the above as follows.
\beD[Stoichiometric and chemical representation]\label{d:crn}
A stoichiometric representation of $f$ is a pair $(\Gamma, r)$ such that
\[
f(x) = \Gamma r(x), \qquad r(x) \ge 0,
\]
where $\Gamma \in \mathbb{R}^{n \times n_R}$ is constant and $r : \mathbb{R}^n_{\ge 0} \to \mathbb{R}^{n_R}_{\ge 0}$ is locally Lipschitz.

It is called chemical if
\[
\Gamma_{i\rho} < 0 \implies r_\rho(x) = x_i \tilde{r}_\rho(x), \qquad \tilde{r}_\rho(x) \ge 0,
\]
and, moreover, $r_\rho$ depends only on variables in $\supp(r_\rho) \subseteq \Gamma^-_i$ whenever $\Gamma_{i\rho} < 0$.
Furthermore,  we assume that each rate $r_\rho(x)$ is monotone non-decreasing with respect to each variable $x_k$ for $k \in \supp(r_\rho)$.
\eeD

\beD[Identification of consumption and production supports]\label{d:Gam}
Let $f$ be chemical in the sense of Definition~\ref{d:crn}. For each $i$, the sets
\bea&&
\Gam_i := \{\, j : f_i^-(x)\ \text{depends nontrivially on } x_j \,\}=\bigcup_{\rho:\,\Gamma_{i\rho}<0}\supp(r_\rho),
\\&&
\Gap_i := \{\, j : f_i^+(x)\ \text{depends nontrivially on } x_j \,\}=\bigcup_{\rho:\,\Gamma_{i\rho}>0}\supp(r_\rho)
\eea
are called consumption and production supports.
\eeD

\beAs[stoichiometric representability]\lbl{as:Gam}
All systems considered in this paper are assumed to
admit a chemical stoichiometric representation (CSR) as above.
\eeAs
   \im We introduce the third ``sources-products" parametrization,   the one most used in  CRNT, by an example. For SIRS,  it is:
   \be{SIRRN}\left(\begin{array}{c|c}
          \text{RN} & \text{rts} \\ \hline
           \text{i}+\text{s}\to 2\text{i} & \beta\, is \\
         \text{i}\to\text{r} & \gamma_i\, i \\
         \text{r}\to\text{s} & \gamma_r\, r \\
         \text{s}\to\text{r} & \gamma_s\, s \\
           \end{array}\right)
\ee

Here each column of $\G$ is    replaced  by a ``reaction" (or ``interaction". For example, the second column in \eqr{SIRc}, representing a transfer from the row of $i$ towards that of $r$,  is replaced
 by $i->r$.

 The $s+i\to 2 i$ in  the first row is the CRNT ``mass-action representation" (not reviewed here). Note  that
 this  provides  a quite natural phenomenological representation of infections at the individual level: the meeting of  a susceptible and  an infectious results in  two infectious.
 \beR
 Note that the (RN,rts) representation of an ODE is not unique. These representations (any of them) have however the advantage  of being the only ones which take into account  the positivity of an ODE (for polynomial ODEs, positivity is equivalent to the existence of a mass-action representation).
 \eeR
\EEN

\subsection{Some interactions between CRN and epidemiology}

CRNT has focused considerably  on graph structures associated with  reactions networks, and on simplifications available when these graphs are (weakly) reversible.
On the other hand, the most encountered infection reaction in epidemiology, $ ``s" + ``i" \to 2``i"$ is irreversible, which explains maybe why
classic CRN and ME
have  not  interacted much in the past (but
  SI, SIS, SIRS models,  appear sometimes in the  CRN literature for the purpose of illustrating various results about autocatalytic models -- see for example \cite{ShiFei,ACK}).

  However, the key fact that Epidemiologic Strains, Minimal Self-Replicable Siphons, and Autocatalytic Cores are conceptually related seems to have passed unnoticed until
  \cite{AABH25}.

We believe that  more  information transfer between CRN and ME could be useful to both, as well as to the other positive system  subfields.
 Let us illustrate this by a few examples:

\BEN \im The {\bf positivity criterion} of \cite{hun} seems unfortunately  unknown outside CRNT, and has been reproved in particular cases an uncountable number of times.

\im Especially relevant is the notion of {\bf siphons}. These  are in one-to-one correspondence with forward invariant boundary faces, cf. Angeli et al.~\cite{AdLS},
and in fact this yields their simplest definition. However, for completeness, we add
(and complete for our purpose) the definition given in the CRN literature:

\beD(siphon, minimal siphon, total/DFE siphon,  critical  siphon, and  DFE  closure/infection)\lbl{d:sif}
\BEN
\im A subset $S$ of species is a siphon/semi-locking set if every reaction that produces a species in S also consumes at least one species in $S$.

\im A siphon is said to be minimal if it does not properly contain any
other siphon.

\im The union of all minimal siphons will be called the total/DFE siphon.

\im A siphon is said to be critical if it does not contain the support of any
linear conservation (``P-flow").
\im The  set of all species which are automatically zero when all the species
in the  DFE  are zero will be called DFE-closure siphon, and any eventual fixed point on the corresponding face will be called DFE.
\EEN
\eeD

 Siphons  are also related to the existence and stability of fixed points on the boundary of the positive orthant.
Indeed, Angeli et al.~\cite{AdLS} proved that if a boundary point is an $\omega$-$limit$ point, then the species set $\Sigma$ in which every species has zero concentration must be a siphon.

\beP \cite[Prop 1]{AdLS}, \cite[Thm 2.5]{AndGAS}, \cite[Lem 2.1]{JohnSie}, \cite[Prop. 3.1]{FelShiu} If an $\omega$-limit $w$ of a strictly positive initial point belongs to a face
$$w \in L_\Sigma:=\{X \in \mathbb{R}^{n}_{\geq 0}:x_i=0, i\in \Sigma\},$$
then $\Sigma$ is a siphon.
\eeP

\beR While the total/DFE siphon is clearly a siphon (since its corresponding face is the intersection of all minimal siphons, hence a siphon), the DFE-closure needs not be a priori  a siphon.  However, it is so in all epidemic models,  because in epidemic models a DFE must  exist, and be an $\omega$-limit  of  strictly positive initial points when $R_0<1$, where $R_0$ is the basic reproduction number defined  as the spectral radius of the \NGM (this siphon is also   critical, since it is unstable when
$R_0>1$, and hence criticality holds by a further result of \cite{AdLS}).\eeR

Our Mathematica package starts any analysis by  inputting the model as a set of reactions, then determines the minimal siphons (using the Julia algorithm provided by Vincent Du \cite{loman2023catalyst}) as a first, easier step towards computing \bfp s residing on each siphon-face \cite{AH26}.

\im
The remarkable NGM theorem of ME states that the  stability  of the \DFE\ (DFE)
 fixed point of the system
 is guaranteed by
 $$R_0 <1.$$

 This result has fascinated ME researchers (a  review \cite{guerra2017basic} focusing solely on measles identified more than 10000 citations that employ $R_0$), due to its probabilistic interpretations, and also because  of the significant simplification \wrt\ the \RH\
 conditions. Despite their  importance, $R_0$ and NGM  seem to have never been mentioned in CRN theory, until \cite{AAN,AAHJ}.  Recently, we have shown  that the natural condition for the NGM theorem to hold is that the boundary face under investigation, the DFE,  is a siphon \cite{AAH26,AH26}.

\EEN

\subsection{Notation and preliminaries}
Let $\mathbb{R}^n$ denote the space of $n \times 1$ column vectors and
$\|\cdot\|_2$ the Euclidean norm. The set of column vectors with positive
entries is denoted by $\mathbb R_+^n$, while the set of column vectors with
nonnegative entries is denoted by $\mathbb R_{\ge 0}^n$. The set of $n\times
m$ matrices with real entries is denoted by $\mathbb R^{n\times m}$.  The
family $\{e_1,\ldots, e_n\}$ denotes the canonical basis of the vector
space $\mathbb{R}^n$. $A^t$ denotes the transpose of $A$, and $A^{-t}$ the
transpose of the inverse of $A$.  If $z \in \mathbb{R}^n$, we denote by $z_i$ the
$i$th component of $z$ (equivalently $z_i = e_i^t . z $).  For a matrix
$A$, we denote by $A(i,j)$ the entry in row $i$, column $j$.  For matrices
$A$ and $B$, we write $A \leq B$ if $A(i,j) \leq B(i,j)$ for all $i$ and
$j$, $A<B$ if $A \leq B$ and $A \neq B$, and $A \ll B$ if $A(i,j) < B(i,j)$
for all $i$ and $j$.

\beD[spectral radius] A) The spectral radius of a matrix $A$ is defined by
$$\rho(A) = max\{|\la| , \la \in Sp(A)\},$$
where $Sp(A)$ denotes the
spectrum of $A$.

B) The largest real part among the eigenvalues of $A$ is denoted by
$$s(A) = \max \left\{{\rm Re}(\lambda)\ : \ \lambda\in Sp(A)\right\}.$$

C) A matrix $A$ is said to be \textit{Hurwitz} if $s(A) < 0$.
\eeD

\beD[Metzler matrix] A Metzler matrix A is a matrix such that $i \neq j \Lra A(i, j) \ge 0.$ These
matrices are also called quasi-positive matrices. 
\eeD

$\vec b, \vec \beta,$ ...  will be used  to denote row vectors. Scalar functions like the logarithm and quotients will be applied to vectors componentwise:
$$\ln y=(\ln y_1,\ln y_2, ...), \fr {1}{y}=\left(\fr 1{y_1},\fr 1{y_2}, \ldots \right),
\fr {y^*}{y}=\left(\fr {y_{1}^*}{y_1},\fr {y_{2}^*}{y_2}, \ldots \right).$$


\medskip
\noindent\textbf{State variables and ambient state space.}
The full state is denoted by
\[
X=(S,I)\in\R_+^m\times\R_+^n,
\]
where
\begin{itemize}
\item $I=(I_1,\ldots,I_n)\in\R_+^n$ denotes the infection variables
(exposed, infected, bacterial, \dots\ compartments),
\item $S=(S_1,\ldots,S_m)\in\R_+^m$ denotes the susceptible-like or resident variables
(susceptible, partially immune, \dots\ compartments).
\end{itemize}
When the state is partitioned relative to a siphon face $\mathcal F_\Sigma$, we use instead
\[
z=(x_\Sigma,y_\Sigma)\in\R_+^N,
\]
where $x_\Sigma$ denotes the variables transversal to $\mathcal F_\Sigma$ and
$y_\Sigma$ the variables tangent to $\mathcal F_\Sigma$.

\subsection{Contents}

In Section \ref{s:ABa}, we revisit the general class of  \ABP s, which include the SAIRS, SEAIR, and SIR-PH models in~\cite{AAH22,AABBGH}. One notable reason for studying these models is that here there exists a natural choice for the matrices from the regular splitting decomposition of \cite{Varga}, \cNGM, specified by the model parameters. These models have a unique interior fixed point, which is also the case for quadratic \gLV s. This is a fascinating aspect, which was investigated in papers such as \cite{Fall,IggidrCEP,Bonzi}, which consider  models with a specific algebraic structure, \cite{ShuaiVan,Shuai13}, which propose a graph-theoretic Lyapunov function construction, \cite{Ye}, which proposes a Poincar\'e-Hopf type approach based on identifying an invariant compact manifold, and \cite{BasarPos}, which makes use of the theory of positive systems. Despite all these efforts, it seems fair to say that a full understanding of the boundary between unique and non-unique interior fixed points is still elusive.

In Section \ref{s:7L} we identify, for the class of  \ABa s without $I$ to $S$ linear feedback, seven laws of the DFE -- EE relay   (beyond rank one,  very little is  known currently).


   In Section \ref{s:7LC}, we describe the effect of a feedback $C \x$ on the seven laws.

Section \ref{s:KP}  recall and applies the computation of Perron eigenvectors
 via Kirchhoff's matrix-tree theorem, which extends the rank-one approach of the previous sections.

 Section~\ref{s:con} provides some conclusions and open
questions.

\def\y{S} \def\x{I} \def\ye{\bar{S}}\def\xe{\bar{I}}

\section{Balanced bilinear models} \label{s:ABa}

The works of \cite{Fall,IggidrCEP,Bonzi}, followed  by \cite{Shuai13},  building on the celebrated  \cite{Van},  suggested  that several connections between the stability  of the DFE  and that of the EE  might exist. In retrospect, this connection is not surprising, since
  the DFE and EE are often members of a ``stability relay team" --- one takes over the stability from the other --see \cite{AH26} for further remarks in this direction. This connection translates into a list of seven laws concerning both the EE and the DFE for \ABa s presented below,  which use the same common building bricks.

 Let us mention that some of these laws are continually being reproved in particular examples in the ME literature, without making reference to the common structure of \ABP s \cite{IggidrCEP,Fall,Bonzi}, which shows the interest in the unification we propose.

 The models we study here are of the form:
 \beD[balanced bilinear model]\label{d:ABP}
A positive ODE system is called a \emph{balanced bilinear model} if, after a permutation of coordinates, it can be written as
\begin{equation}\label{ABP}
\begin{cases}
\y'= \f(\y)-\diag(\y)\,\bbe \x + C\x, \; \f(\y)=\Lambda +A_S \y\\[2mm]
\x'=P\,\diag(\y)\,\bbe \x + A\x,
\end{cases}
\end{equation}
where
\[
\y\in\R_+^m,\qquad \x\in\R_+^n,
\]
and:
\begin{enumerate}
\item $I \in \mathbb{R}^n_+$ are $n$ DFE closure/infection compartments.
 \item $S \in \mathbb{R}^m_+$ are $m$ non-infection/susceptible compartments (the complement of $\x$).
 \item  \(C\in\R_+^{m\times n}\) yields feedback from $\x$ to $\y$.
 \item $\f(\y) $ is the dynamics in the $\y$ \com s, in the absence of the $\x$ \com s. Below we always assume  $\f(\y)=\Lambda +A_S \y, \La\in \mathbb{R}_+^m,  A_S \in \mathbb{R}^{m,m}$.

 \item $A_S \in \mathbb{R}^{m\times m}, A\in \mathbb{R}^{n,n}$ are stable Metzler matrices modeling the intra flows between S and I compartments, \resp. The second is precisely the matrix appearing in the NGM formula $K=F(-A)^{-1}$.
 \item $\bbe \in \mathbb{R}_+^{m\times n} > 0$ is the WAIFW (Who Acquires Infection From Whom) matrix, and $ \bbe (i,j) $ denotes the contact and infectivity rate of $ \x _j $ within the group $ \y_i $. This matrix, which exists only for models which are bilinear in $\y,\x$, is related to the square matrix $F$ from the NGM formula.

 \item $P\in \mathbb{R}_+^{n\times m}$ is a column-stochastic matrix that distributes the vector $\diag(\y)\,\bbe \x$ among the $\x$ \com s of the infection.  The fact that $1 .P=1$ implies that the bilinear terms cancel out in the sum of all equations, suggesting the name ``balanced bilinear" for this class.

\end{enumerate}

\eeD

\begin{hypothesis} \label{hyp: bilinear_model}
\begin{itemize}
\item Any susceptible compartment is accessible from a susceptible compartment with inflow/recruitment.

\item Any infected-infectious compartment is accessible from at least one entry-point for infection.
\end{itemize}
\end{hypothesis}

\beD[linear forces of infection and bilinear rates of infection] Below, a function \be{F}F(\y):=P\,\diag(\y)\,\bbe=\sum_{i=1}^m \y_i\,p_i\,\vec b_i \quad  \in  \mathbb{R}_+^{n\times n}, \quad \y \in \mathbb{R}^m_{+},\ee where $P, \bbe$ satisfy the conditions of definition \ref{d:ABP}, and $p_i$, $\vec b_i$ are the $i$-th column of $P$ and   $i$-th row of $\bbe$, \resp, will  be called a
linear  force of infection, and
$$F(\y) \, \x=P\,\diag(\y)\,\bbe \x =P  \diag (\bbe \x )\y\quad  \in \mathbb{R}^n_{+}$$
will be called a
bilinear  rate of infection.

 \eeD

 \beR[some alternative representations for  bilinear rates of infection $F(\y)$]\label{rem:Falts}

A bilinear rate of infection \mbw\
 \be{altF}F(\y)\,\x=P\,\diag(\y)\,\bbe \x=P\,(\y\,\odot\bbe \x)=P  \diag (\bbe \x )\y,\ee

 where $\odot$ denotes the componentwise (Hadamard) product, and the identities
   follow from the standard vector identity $\diag(a)\,b = a\odot b = \diag(b)\,a$.

   The {\bf ``Diag switch"} equality between the second and forth will be used below.

\eeR

The following two formulas hold independently of the dimension of $\y$ and of rank one assumptions:
\BEN \im  If $A_S$ is Metzler stable (ensuring the element-wise positivity of $(-A_S)^{-1})$), then
the DFE is unique, given by:
$$\y_0=(-A_S)^{-1}\Lambda\in \mathbb{R}_+^m.$$

\im
the basic reproduction number is
$${R}_0 = \rho\left(K(\y_0)\right)= \rho\left(F(\y_0) \pr{-A}^{-1}\right),$$
where $\y_0$ is the disease-free equilibrium.
\EEN

\ssec{Rank one models of P-type and B-type}

 \beR  \lbl{r:r1}
 Several of the
explicit formulas in Section \ref{s:7L}  hold only when either $P=\ba_n\, (1,1, ...,1)$ or $\bbe=\ba_m\, \vbe$ are of rank one.
   Since the coefficients $\y_i$ vary independently in $\mathbb{R}^m_{+}$,
    it follows that $\mathrm{rank}(F(\y))=1$ for all interior $\y$
   if and only if either $\mathrm{rank}(P)=1$ or $\mathrm{rank}(\bbe)=1$.
We have thus two (non-disjoint) cases:

\begin{itemize}
 \item \emph{Case~(P):} $P=\ba_n\,\vec{\bff 1}$,\quad $\ba_n\in\mathbb{R}^n_+$,\ $\vec{\bff 1}\ba_n=1$;\quad $F(\y)\x=\ba_n\,(\vec\y\bbe\,\x)$.
 \item \emph{Case~(B):} $\bbe=\ba_m\vbe$,\quad $\ba_m\in\mathbb{R}^m_+$,$\vec{\bff 1}\ba_m=1$,\ $\vbe\in\mathbb{R}^{1\times n}_+$;\quad $F(\y)\x=P\diag(\y)\ba_m\cdot(\vbe\,\x)$.
\end{itemize}

In both cases the \ABP\  \eqref{ABP} without feedback from $\x$ to $\y$   reads
\be{FIS0}
\bc
\y'=\f(\y)-\diag[\y]\,\bbe\x+ C\x, 
\\[2pt]
\x'=\pp{F(\y)+A}\x:=\pp{K(\y)-Id_n} (-A) \x,\qquad K(\y):=F(\y)(-A)^{-1},
\ec
\ee
where $K(\y)$ is the NGM.
\eeR
In Case~(P) once infected, individuals enter the infectious
compartments according to the stochastic vector $\ba$, regardless of their susceptible class
of origin. In Case~(B), the infectivity row $\vbe$ is common to all susceptible classes,
which differ only through $\ba_m$.

\beR[Rank-one factorizations of $K$ and $\T K$]\label{r:NGMeig}
Recall
\(
F(\y)=P\,\diag(\y)\,\bbe ,
\)
and
\[
\T K(\y_0)=(-A)^{-1}F(\y_0),\qquad
K(\y_0)=F(\y_0)(-A)^{-1}=(-A)\,\T K(\y_0)\,(-A)^{-1},
\]
so $K(\y_0)$ and $\T K(\y_0)$ are similar.

\medskip

\noindent
\textbf{Case~(P).}
$
P=\ba_n\,\vec{\bff 1},
\qquad
\vec{\bff 1}\ba_n=1,$ implying
$ F(\y_0)=\ba_n\,\y_0\,\bbe,
$
and therefore
\[
K(\y_0)=\ba_n\,\vec\pi_K,
\qquad
\vec\pi_K:=\y_0\,\bbe(-A)^{-1}.
\]
Hence one may choose Perron eigenvectors of $K(\y_0)$ as
\[
\text{right: } \mathbf w_K=\ba_n,
\qquad
\text{left: } \vec\pi_K=\y_0\,\bbe(-A)^{-1}.
\]
For $\T K(\y_0)$ one has
\[
\T K(\y_0)=\mathbf w_{\T K}\,\vec\pi_{\T K},
\qquad
\text{right: }\; \mathbf w_{\T K}:=(-A)^{-1}\ba_n,
\qquad
\text{left: }\; \vec\pi_{\T K}:=\y_0\,\bbe.
\]

\medskip

\noindent
\textbf{Case~(B).} $
\bbe=\ba_m\vbe,
\qquad
\vec{\bff 1}\ba_m=1,
\qquad
\vbe\in\mathbb R^{1\times n}_+, $ implying
$F(\y_0)=\mathbf w_F\,\vbe,
\qquad
\mathbf w_F:=P\,\diag(\y_0)\,\ba_m,
$
and therefore
\[
K(\y_0)=\mathbf w_K\,\vec\pi_K,
\qquad
\mathbf w_K:=P\,\diag(\y_0)\,\ba_m,
\qquad
\vec\pi_K:=\vbe(-A)^{-1}.
\]

Summarizing,  one may choose Perron eigenvectors as follows:
 \begin{table}[H]
  \centering
  \renewcommand{\arraystretch}{2.0}
  \begin{tabular}{llllll}
  \hline
  Case & Assumption
    & \multicolumn{2}{c}{$K(\y_0)$}
    & \multicolumn{2}{c}{$\T K(\y_0)=\mathbf{w}_{\T K}\,\vec\pi_{\T K}$} \\
  \cline{3-4}\cline{5-6}
  & & right $\mathbf{w}_K$ & left $\vec\pi_K$
    & right $\mathbf{w}_{\T K}$ & left $\vec\pi_{\T K}$ \\
  \hline
  (P) & $P=\ba_n\,\vec{\bff 1}$
    & $\ba_n$
    & $\y_0\,\bbe(-A)^{-1}$
    & $(-A)^{-1}\ba_n$
    & $\y_0\,\bbe$ \\
  (B) & $\bbe=\ba_m\vbe$
    & $P\,\diag(\y_0)\,\ba_m$
    & $\vbe(-A)^{-1}$
    & $(-A)^{-1}P\,\diag(\y_0)\,\ba_m$
    & $\vbe$ \\
  \hline
  \end{tabular}
  \end{table}
\medskip

\noindent
In both cases,
\[
R_0=\rho\bigl(K(\y_0)\bigr)=\rho\bigl(\T K(\y_0)\bigr)
=\vec\pi_K\,\mathbf w_K
=\vec\pi_{\T K}\,\mathbf w_{\T K}.
\]
Moreover, the right and left Perron eigenvectors of $K(\y_0)$ and $\T K(\y_0)$ are related by
\[
\mathbf w_K=(-A)\,\mathbf w_{\T K},
\qquad
\vec\pi_{\T K}=\vec\pi_K\,(-A).
\]
\eeR

\subsection{Independence of the two rank--one assumptions:ExP.wl,ExB.wl}

Recall that the infection dynamics in \eqref{FIS0} may be written
\[
\x' = F(\y)\x + A\x,
\qquad
F(\y):=P\,\diag(\y)\,\bbe ,
\]
with simplifying cases
\[
\text{(P)}\quad P=\ba_n(\vec{\bff 1}),
\qquad
\text{(B)}\quad \bbe=\ba_m\,\vbe.
\]
We show now that the two assumptions are algebraically independent.

\begin{lemma}[Rigidity of the factorization]
\label{lem:rigidity}
Let $F(\y)=P\diag(\y)\bbe$.
Assume that every row of $\bbe$ and of $\widetilde{\bbe}$ is nonzero.
If two column-stochastic pairs $(P,\bbe)$ and $(\widetilde{P},\widetilde{\bbe})$
generate the same infection operator for all $\y$,
\[
P\diag(\y)\bbe=\widetilde{P}\diag(\y)\widetilde{\bbe}
\qquad\forall\,\y\in\mathbb R^m,
\]
then necessarily
\[
P=\widetilde{P},\qquad \bbe=\widetilde{\bbe}.
\]
\end{lemma}

\begin{proof}
Write
\[
P\diag(\y)\bbe
=\sum_{i=1}^m \y_i\,p_i\,\vec b_i,
\]
where $p_i$ is the $i$th column of $P$ and $\vec b_i$ the $i$th row of $\bbe$.
The same expansion holds for $(\widetilde{P},\widetilde{\bbe})$:
\be{str}
\sum_{i=1}^m \y_i\,p_i\,\vec b_i
=
\sum_{i=1}^m \y_i\,\tilde p_i\,{\tilde \vec b_i}
\qquad\forall\,\y.
\ee
Since the coefficients $\y_i$ are independent, we obtain
\[
p_i\,\vec b_i=\tilde p_i\,{\tilde \vec b_i}
\qquad\text{for each }i.
\]
Since $P$ and $\widetilde{P}$ are column-stochastic, each column satisfies
$\vec{\bff 1}p_i=1$ and $\vec{\bff 1}\tilde p_i=1$.
Because $\vec b_i\neq0$ and $\tilde{\vec b_i}\neq0$, the outer products above have rank one, hence
\[
p_i=\lambda_i\tilde p_i
\]
for some scalar $\lambda_i>0$.
Applying both unit-sum conditions gives $\lambda_i=1$, hence $p_i=\tilde p_i$.
It then follows from $p_i(\vec b_i-{\tilde \vec b_i})=0$ and $p_i\neq 0$ that
$\vec b_i={\tilde \vec b_i}$.
Hence $P=\widetilde{P}$ and $\bbe=\widetilde{\bbe}$.
\end{proof}

This lemma shows that the infection dynamics uniquely determines the
column-stochastic pair $(P,\bbe)$, so one cannot transform a model satisfying one of the
two assumptions into the other without changing the dynamics.

\beT[Existence of Case (P) without Case (B)]
\label{prop:P-not-beta}
There exist systems of the form \eqref{ABP} with
$P=\ba_n(\vec{\bff 1})$ but $\bbe$ not of rank one.
\eeT

\begin{proof}
Consider $m=n=2$ and
\[
\ba_n=
\begin{pmatrix}
1/3\\
2/3
\end{pmatrix},
\qquad
P=\ba_n(1,1)=
\begin{pmatrix}
1/3&1/3\\
2/3&2/3
\end{pmatrix}.
\]
Let
\[
\bbe=
\begin{pmatrix}
1&2\\
3&4
\end{pmatrix}.
\]
Then
\[
\det\bbe=-2\neq0,
\]
so $\mathrm{rank}(\bbe)=2$.
Hence the model satisfies Case~(P) but not Case~(B).

By Lemma \ref{lem:rigidity}, no alternative rank-one matrix
$\widetilde{\bbe}$ can reproduce the same infection operator
$F(\y)=P\diag(\y)\bbe$.
\end{proof}

\beT[Existence of Case (B) without Case (P)]
\label{prop:beta-not-P}
There exist systems of the form \eqref{ABP} with $\bbe$ rank one but
with no representation satisfying $P=\ba_n(\vec{\bff 1})$.
\eeT

\begin{proof}
Again take $m=n=2$ and define
\[
\ba_m=
\begin{pmatrix}
1/3\\
2/3
\end{pmatrix},
\qquad
\vbe=\begin{pmatrix}3&9\end{pmatrix},
\qquad
\bbe=\ba_m\,\vbe=
\begin{pmatrix}
1&3\\
2&6
\end{pmatrix}.
\]
Thus $\bbe$ has rank one.

Now take
\[
P=
\begin{pmatrix}
1&0\\
0&1
\end{pmatrix}.
\]
This matrix is column-stochastic but its columns are not proportional,
so it cannot be written in the form $\ba_n(\vec{\bff 1})$.

Suppose that the same infection operator could be written using
$\widetilde{P}=\ba_n(\vec{\bff 1})$.
Then since B is common to both representations; 
\[
P\diag(\y)\bbe=\widetilde{P}\diag(\y)\bbe
\qquad\forall\,\y.
\]
By Lemma \ref{lem:rigidity} this would imply $P=\widetilde{P}$, which is
impossible because the columns of $P$ are not proportional.

Hence Case~(B) does not imply Case~(P).
\end{proof}

{\bf Conclusion}:
The two simplifying assumptions
\[
P=\ba_n(\vec{\bff 1}),
\qquad
\bbe=\ba_m\,\vbe
\]
are independent structural properties of the bilinear epidemic model.
Neither can be deduced from the other without changing the infection
dynamics.


\ssec{Perron-Frobenius strong threshold theorem for ABP models without feedback}

 The following result, motivated by \ABP s introduced in Section \ref{s:ABa},  is independent of the dimension of $\y$ and does not rely on rank-one assumptions. Previous constructions of endemic equilibria in this spirit were formulated only for particular cases; see for example \cite{Bonzi,IggidrCEP,Shuai13}, where explicit rank-one representations are used.

We show that such constructions may be replaced by a two-step procedure based on Perron--Frobenius theory. First, one determines $S$ by solving the scalar equation
\[
\rho(\widetilde K(S))=1,
\]
where $\widetilde K(S)$ is the next-generation operator associated with the infected subsystem. Second, one recovers the infected equilibrium direction from a positive left Perron eigenvector of $\widetilde K(S)$, and determines its amplitude through the susceptible subsystem.

This yields a constructive characterization of endemic equilibria under the sole assumption that the infected Jacobian block admits a regular splitting and is irreducible. While the argument extends to all epidemiological models satisfying these conditions, we formulate it here for \ABP s in order to relate it to subsequent results in the rank-one setting.
 \beT[Perron-Frobenius  strong threshold theorem for ABP models without feedback]\label{thm:TK_DFE}
Fix \(\y\in\R^m_{\ge 0}\), and consider the linear equilibrium equation
\begin{equation}\label{eq:linx_again}
0=A\x+P\,\diag(\y)\,\bbe\x,
\end{equation}
where \(A\) is Metzler--Hurwitz and \(P,\bbe\ge 0\). Introduce:
\[
\T K(\y):=\widetilde G\,\diag(\y),
\qquad
\widetilde G:=\bbe(-A)^{-1}P,
\]
and assume that \(\widetilde G\) is irreducible. Then:

\begin{enumerate}
\item[(i)] {For the fixed $\y$,} Equation \eqref{eq:linx_again} admits a nonzero solution \(\x\blue{>} 0\) if and only if 
\[
\rho\big(\T K(\y)\big)=1.
\]
{In particular, \(\x=0\) is always a solution of \eqref{eq:linx_again}. Moreover, it is the unique nonnegative solution  if and only if $\rho\big(\T K(\y)\big)< 1$.  }

\item[(ii)] If \(\rho(\T K(\y))=1\) and \(\y\gg0\), then all nonzero nonnegative solutions of \eqref{eq:linx_again} form a one-dimensional positive ray:
\[
\x=k \,(-A)^{-1}P\,\diag(\y)\,v(\y),\qquad k>0,
\]
where \(v(\y)\gg0\) is a Perron eigenvector of \(\T K(\y)\), unique up to multiplication by a positive scalar.

\item[(iii)] Let \(\y_0\gg0\) be the disease-free value, and set
\[
R_0:=\rho\big(\T K(\y_0)\big).
\]
If \(R_0<1\), then for every \(0\le \y\le \y_0\),
\[
\rho\big(\T K(\y)\big)\le R_0<1.
\]
Hence, for every \(0\le \y\le \y_0\), equation \eqref{eq:linx_again} admits only the trivial nonnegative solution \(\x=0\).

\item[(iv)] If \(R_0=1\), then \(\y_0\) itself lies on the threshold hypersurface
\[
\rho\big(\T K(\y)\big)=1.
\]

\item[(v)] If \(R_0>1\), then there exists a unique equilibrium point on the ray \(t\y_0\), \(t\ge 0\), with 
\[
\bar\y:=\frac{1}{R_0}\,\y_0,
\]
and moreover
\[
\rho\big(\T K(\bar\y)\big)=1.
\]

\item[(vi)] More generally, if \(\y^*\gg0\) satisfies
\[
\rho\big(\T K(\y^*)\big)=1,
\]
then the \(\x\)-equation determines a unique positive ray
\[
\x=k\,\x_*(\y^*),\qquad k>0,
\]
where
\[
\x_*(\y^*):=(-A)^{-1}P\,\diag(\y^*)\,v(\y^*).
\]
\end{enumerate}
\eeT

\beR[$\T K$ is  cyclic permutation of the next-generation matrix]\label{r:cyclic}
At \(\y=\y_0\), the operator \(\T K(\y_0)\) is the cyclically permuted form of the
usual next-generation matrix
\[
K(\y_0)=P\,\diag(\y_0)\,\bbe(-A)^{-1}.
\]
Indeed,
\[
\T K(\y_0)=\bbe(-A)^{-1}P\,\diag(\y_0),
\]
and the matrices \(K(\y_0)\) and \(\T K(\y_0)\) have the same nonzero spectrum. Hence
\[
R_0=\rho\big(K(\y_0)\big)=\rho\big(\T K(\y_0)\big).
\]
This is why Theorem~\ref{thm:TK_DFE} applies also with the standard \NGM\ replacing $\T K$.
\eeR

\beR[Meaning of the threshold ray]\label{rem:threshold-ray}
Fixing an  \(\y^*\gg0\) such that
\[\bc
0=A\x+P\,\diag(\y^*)\,\bbe\x\\
\rho(\T K(\y^*))=1
\ec
\]
like in the last subpoint does not determine a unique equilibrium vector \(\x\), but only a one-dimensional
positive ray
\[
\x = k\,\x_*(\y^*),\qquad k>0,
\]
where \(\x_*(\y^*)\) is fixed up to normalization.

In particular:
\begin{itemize}
\item the \emph{profile} of the infected state (i.e.\ the ratios between components of \(\x\))
      is uniquely determined by the Perron eigenvector of \(\T K(\y^*)\);
\item the \emph{magnitude} of \(\x\) is not determined at the linear level, and remains
      free up to a positive scalar factor.
\end{itemize}

Thus, at threshold, the system selects a unique infection pattern but not its amplitude.
The latter is fixed only when coupling with the nonlinear \(\y\)-equation.
\end{remark}

\begin{proof}
Since \(A\) is Hurwitz, it is invertible, and since \(A\) is Metzler, one has
\[
(-A)^{-1}\ge 0.
\]
Thus \eqref{eq:linx_again} is equivalent to
\begin{equation}\label{eq:xfixed_new}
\x=(-A)^{-1}P\,\diag(\y)\,\bbe\x.
\end{equation}
Applying \(\bbe\) gives
\begin{equation}\label{eq:ufixed_new}
u=\T K(\y)\,u,\qquad u:=\bbe\x\ge 0.
\end{equation}

To prove (i), suppose first that \eqref{eq:linx_again} admits a nonzero solution \(\x\ge 0\). Then \(u=\bbe\x\neq0\), and \eqref{eq:ufixed_new} shows that \(1\) is an eigenvalue of the nonnegative matrix \(\T K(\y)\) with a nonzero nonnegative eigenvector. Hence
\[
\rho\big(\T K(\y)\big)=1.
\]
Conversely, assume that \(\rho(\T K(\y))=1\), and let \(v(\y)\ge 0\), \(v(\y)\neq0\), satisfy
\[
\T K(\y)\,v(\y)=v(\y).
\]
Define
\[
\x^*:=(-A)^{-1}P\,\diag(\y)\,v(\y)\ge 0.
\]
Then
\[
\bbe\x^*=\T K(\y)\,v(\y)=v(\y),
\]
and therefore \eqref{eq:xfixed_new} holds for \(\x^*\). Thus \(\x^*\) is a nonzero nonnegative solution of \eqref{eq:linx_again}. This proves (i).

For (ii), assume that \(\rho(\T K(\y))=1\) and \(\y\gg0\). Since
\[
\T K(\y)=\widetilde G\,\diag(\y),
\]
and \(\diag(\y)\) is positive diagonal, \(\T K(\y)\) has the same zero pattern as \(\widetilde G\), hence is irreducible. By Perron--Frobenius, the eigenspace of \(\T K(\y)\) for the eigenvalue \(1=\rho(\T K(\y))\) is one-dimensional and generated by a vector \(v(\y)\gg0\), unique up to multiplication by a positive scalar.

Let \(\x\ge 0\), \(\x\neq0\), solve \eqref{eq:linx_again}, and set \(u=\bbe\x\). Then \(u\neq0\) and \eqref{eq:ufixed_new} gives
\[
u=\T K(\y)\,u.
\]
Hence \(u=k\,v(\y)\) for some \(k>0\). Substituting into \eqref{eq:xfixed_new} yields
\[
\x=k\,(-A)^{-1}P\,\diag(\y)\,v(\y).
\]
Conversely, every vector of this form satisfies \eqref{eq:xfixed_new}. This proves (ii).

For (iii), if \(0\le\y\le\y_0\), then
\[
\diag(\y)\le\diag(\y_0),
\]
hence
\[
\T K(\y)=\widetilde G\,\diag(\y)\le \widetilde G\,\diag(\y_0)=\T K(\y_0)
\]
entrywise, because \(\widetilde G\ge 0\). By monotonicity of the spectral radius on nonnegative matrices,
\[
\rho\big(\T K(\y)\big)\le \rho\big(\T K(\y_0)\big)=R_0<1.
\]
The last statement follows from (i).

Part (iv) is immediate from the definition of \(R_0\).

For (v), if \(t\ge 0\), then
\[
\T K(t\y_0)=\widetilde G\,\diag(t\y_0)=t\,\widetilde G\,\diag(\y_0)=t\,\T K(\y_0),
\]
so
\[
\rho\big(\T K(t\y_0)\big)=t\,R_0.
\]
Hence \(\rho(\T K(t\y_0))=1\) if and only if \(t=1/R_0\), and therefore the threshold hypersurface is met on the ray \(t\y_0\) at the unique point \(\bar\y=\y_0/R_0\).

Finally, (vi) is an immediate application of (ii) with \(\y=\y^*\).
\end{proof}

\beR[DFE corollary for the coupled model]\label{r:DFEcor}
For the coupled balanced bilinear system
\[
\begin{cases}
\y'=\Lambda+A_S\y-\diag(\y)\,\bbe\x+C\x,\\[1mm]
\x'=A\x+P\,\diag(\y)\,\bbe\x,
\end{cases}
\qquad
\y_0=-A_S^{-1}\Lambda,
\]
with \(A_S\) Metzler--Hurwitz, Theorem~\ref{thm:TK_DFE}(iii) provides the
existence/uniqueness part of the disease-free analysis: if \(R_0<1\), then for every
\(0\le \y\le \y_0\) the \(\x\)-equation admits only the trivial nonnegative solution
\(\x=0\). Thus \((\y_0,0)\) is the unique disease-free equilibrium. The additional
stability statement for the full coupled system is then obtained by linearization of the
\((\y,\x)\)-system at \((\y_0,0)\).
\eeR

\beR[Spectral mechanism and CEP reduction]\label{r:CEP}
Theorem~\ref{thm:TK_DFE} isolates the essential structure of the equilibrium problem:
\[
\rho(\T K(\y))<1 \iff \x=0,
\qquad
\rho(\T K(\y))=1 \iff \exists\,\x\neq0.
\]
Hence every nontrivial equilibrium must lie on the hypersurface
\[
\rho\big(\T K(\y)\big)=1.
\]
On this hypersurface, the \(\x\)-equation determines a unique positive ray
\[
\x=k\,\x_*(\y),\qquad k>0.
\]
Thus the equilibrium problem separates into two steps:
\[
\text{(i) solve } \rho\big(\T K(\y)\big)=1,
\qquad
\text{(ii) determine the scalar }k\text{ from the }\y\text{-equation.}
\]
In particular, if the spectral equation admits a unique solution \(\y^*\gg0\), then the
only remaining unknown is the scalar amplitude \(k\). This is the precise form of the
CEP mechanism beyond rank one.
\eeR

\beR[Loop interpretation]\label{r:loop}
The operator
\[
\T K(\y)=\bbe(-A)^{-1}P\,\diag(\y)
\]
represents one full circulation through the positive feedback loop
\[
\x \;\xrightarrow{\bbe}\; \R^m
\;\xrightarrow{\diag(\y)}\; \R^m
\;\xrightarrow{P}\; \R^n
\;\xrightarrow{(-A)^{-1}}\; \R^n.
\]
If \(z=\bbe\x\), then the equilibrium condition becomes
\[
z=\T K(\y)\,z.
\]
Thus \(\rho(\T K(\y))=1\) expresses exact balance between dissipation through \(A\) and
amplification through the feedback loop. This is the structural origin of the threshold.
\eeR 
\section{Seven explicit formulas for  rank one \ABP s without linear feedback matrix $C$} \label{s:7L}

\subsection{The computation of the DFE and EE for  rank one \ABP s without linear feedback matrix $C$}


\beT[Explicit consequences of the Perron-Frobenius  strong threshold theorem  for rank-one]\label{t:end}

For any rank-one \ABP\ \eqref{FIS0} satisfying Hypothesis~\ref{hyp: bilinear_model},
there exists a unique endemic equilibrium $(\xe,\ye)$ with strictly positive coordinates
if and only if $R_0>1$. Moreover, the following four laws hold:

\BEN
\im
\textbf{First Law: Basic Reproduction Number \cite[(3.2)]{Fall}, \cite[(6)]{Bonzi}}
$R_0$ \mbw:
\be{R0ABa}
R_0=\y_0\mR,
\ee
where $\mR\in\mathbb{R}^m_+$ is the replacement vector of Definition~\ref{def:mR} below.

\beD[Replacement vector]\label{def:mR}
The {\bf replacement vector} $\mR\in\mathbb{R}^m_{\ge 0}$ is given by
\be{mR}
\mR:=\begin{cases}\bbe(-A)^{-1}\ba & \text{Case~(P)},\\[4pt]
\diag(\ba_m)\,[\vbe(-A)^{-1}P]^\top & \text{Case~(B),}\end{cases}
\ee
equivalently in Case~(B): $\mR_i=(\ba_m)_i\,\vbe(-A)^{-1}p_i$, $i=1,\dots,m$,
where $p_i$ is the $i$-th column of $P$.
\eeD

\beR
{$\mR$ is the vector whose $i$th component is the expected number of infections generated by a single susceptible individual in class $i$.}
\eeR

\im
\textbf{Second Law: $\xe$ is proportional to right eigenvector and ``dwell times" \cite[(3.5)]{Fall}}
\be{xee}
\xe = k\,(-A)^{-1}\ba_{\mathrm{eff}},\qquad
\ba_{\mathrm{eff}}:=\begin{cases}\ba_n & \text{Case~(P)},\\ P\diag(\ye)\ba_m & \text{Case~(B).}\end{cases}
\ee
In Case~(P), $D_w:=(-A)^{-1}\ba_n$ is a \emph{fixed} vector (independent of the equilibrium),
so the direction of $\xe$ is constant, determined by these parameters alone.

\begin{definition}[dwell times]
For any \MH\ matrix $A$, or, equivalently, for any finite-space irreducible CTMC with absorbing states and transitions between the
transient states encapsulated in a matrix $A$, dwell times are the (positive) elements of
$(-A)^{-1}$.
\end{definition}
(This name, used in the theory of Markov processes, and which has also been useful in mathematical epidemiology--  see for example \cite{Hurtado,AABJ,AAB}, is explained by the well-known fact that these are precisely the means (averages) of
times before absorption, conditioned on both a starting point and a last point before
absorption).

$D_w=(-A)^{-1}\ba_n$ may be interpreted
as ``expected infection times until recovery'', starting in the first state with law $\ba_n$.
\im
\textbf{Third Law: Endemic $\ye$ as function  of the normalization constant}
\be{yee}
\ye=\pr{\diag[\bbe\xe]-A_S}^{-1}\La
=\pr{k\diag[\bff{a}]-A_S}^{-1}\La,\qquad
\bff{a}:=\begin{cases}\mR & \text{Case~(P)},\\ \ba_m & \text{Case~(B)}\end{cases}
\ee
holds, whenever the intervening matrix is invertible.

\im
\textbf{Fourth Law: $\ye$ normalization \cite[(3.6)]{Fall}, \cite[(4)]{Bonzi}}
Define now, whenever the intervening matrix is invertible, the function
\be{H}
H(k):=\mR\pr{k\diag[\bff{a}]-A_S}^{-1}\La.
\ee
Then
\be{rH}
\rho\!\big(K(\ye(k))\big)=H(k),
\ee
and the normalization constant $k=\ye\bbe\xe$ is the unique positive solution of
\be{mRv}
H(k)=1;
\ee
equivalently,
\[
\ye(k)\mR=1.
\]

When $A_S=-\diag(\mu_S)$  for example,
$H(k)=\displaystyle\sum_j\frac{\mR_j\,\La_j}{k\,a_j+\mu_j}$,
and for $n=m=1$:
\be{n1}
\ye=\frac{1}{\mR},\qquad k=\mu_S\ye(R_0-1)=\La-\frac{\mu_S}{\mR}.
\ee
\EEN

\eeT

\beR
When $m>1$, $k$ is in general not rational, making an explicit determination of
$(\xe,\ye)$ difficult.
\eeR

\begin{proof}
\noindent\textbf{Step 1.}
The argument follows \cite[(3.4--3.7)]{Fall}; we give it uniformly for both cases.
Let $F(\ye)=\bff u_e\,\vec{\bff w_e}$ be the rank-one factorisation at the EE, where

\smallskip\noindent
Case~(P): $\bff u_e=\ba$,\ $\vec{\bff w_e}=\ye\bbe$;\qquad
Case~(B): $\bff u_e=P\diag(\ye)\ba_m$,\ $\vec{\bff w_e}=\vbe$.

\smallskip
\noindent\textbf{Step 2.}
At the EE, $\x'=0$ gives
\[
(F(\ye)+A)\xe=0,
\]
hence
\[
\bff u_e(\vec{\bff w_e}\xe)=-A\xe.
\]
Therefore
\be{xsb}
\xe=k\,(-A)^{-1}\bff u_e,\qquad k:=\vec{\bff w_e}\xe>0,
\ee
which is \eqref{xee}.

\noindent\textbf{Step 3.}
The equation $\y'=0$ gives
\[
\ye=(\diag[\bbe\xe]-A_S)^{-1}\La.
\]
Using
\[
\bbe\xe=k\,\bbe(-A)^{-1}\ba=k\mR \qquad\text{in Case~(P)},
\]
and
\[
\bbe\xe=\ba_m\vbe\xe=k\ba_m \qquad\text{in Case~(B)},
\]
yields
\[
\ye(k)=\pr{k\diag[\bff a]-A_S}^{-1}\La,
\]
which is \eqref{yee}.

\noindent\textbf{Step 4.}
Multiplying \eqref{xsb} by $\vec{\bff w_e}$ gives
\[
k=k\,\vec{\bff w_e}(-A)^{-1}\bff u_e.
\]
Since $k>0$, division by $k$ yields
\[
\vec{\bff w_e}(-A)^{-1}\bff u_e=1.
\]
In both cases this is exactly
\[
\ye\,\mR=1.
\]
Substituting \(\ye=\ye(k)\) gives
\[
H(k)=\mR\pr{k\diag[\bff a]-A_S}^{-1}\La=1.
\]

\noindent\textbf{Step 5.}
Since $K(\ye)$ has rank one, its unique nonzero eigenvalue equals its trace.
In both cases this gives
\[
\rho\!\big(K(\ye)\big)=\ye\,\mR.
\]
Therefore, for \(\ye=\ye(k)\),
\[
\rho\!\big(K(\ye(k))\big)=\vec\mR\pr{k\diag[\bff a]-A_S}^{-1}\La=H(k),
\]
which is \eqref{rH}. Hence the Perron condition
\[
\rho\!\big(K(\ye(k))\big)=1
\]
is equivalent to the normalization law
\[
H(k)=1,
\]
or equivalently \(\ye(k)\mR=1\).

\noindent\textbf{Step 6.}
Uniqueness of the positive solution follows since $H$ is strictly decreasing,
\[
\lim_{k\to\infty}H(k)=0,
\qquad
H(0)=R_0>1.
\]
\end{proof}

\subsection{Fifth Law: An explicit \Lf\ establishing global asymptotic stability of the DFE for
\ABa\ (Case~P) when $R_0<1$}\label{s:LDFE}

The following result is for Case~(P) ($P=\ba_n\vec{\bff 1}$).
\beT[Stability of DFE for (P) \ABa s]
\label{t:DFE} For an $\x$-$\y$ ODE
\bea
\bc
\y'=\f(\y)-\diag[\y]\,\bbe\x,\\
\x'=\pp{\ba_n\y\bbe+A}\x,
\ec
\eea
let $R_0$ be defined by the natural regular splitting of the \jin\ into $\ba_n\y\bbe$ and $A$.
The following holds:
\begin{itemize}
\item If ${R}_0>1$, then the disease-free equilibrium is unstable.

\item If ${R}_0\leq1$, and $A_S=-\diag(\mu_S)$ is a strictly negative diagonal matrix, then
\be{LDFE}
V_{DF}(\x,\y)=R_0\pr{{\vec{\bff 1}}\y-\y_0\ln(\y)}+\y_0\bbe(-A)^{-1}\x
\ee
is a \Lf, whose minimum is attained uniquely at the DFE, which is therefore
globally asymptotically stable.
\end{itemize}
\eeT

\beR[Lyapunov weights  are related to the left NGM eigenvector]\label{r:LDFEPerron}
The function \eqref{LDFE} is of the form
\[
\mu H(\y)+ L(\x), \qquad \mu=R_0,
\]
with tangential Goh--Volterra--Horn--Jackson term
\[
H(\y)=\vec{\bff 1}\y-\y_0\ln(\y)
\]
and
with transversal Perron term
\[
L(\x)=\y_0 \bbe(-A)^{-1},\x,
\]
whose weights
are a left Perron eigenvector of $K(\y_0)=F(\y_0)(-A)^{-1}$ in both cases (P) and (B).
\eeR

\begin{proof}[Proof of Theorem~\ref{t:DFE}]
At the DFE one has
\[
\y_0=(-A_S)^{-1}\Lambda,
\qquad
A_S=-\diag(\mu_S),
\]
hence
\[
\Lambda=\diag(\mu_S)\y_0,
\qquad
\f(\y)=\Lambda+A_S\y=-\diag(\mu_S)(\y-\y_0).
\]

Set
\[
H(\y):={\vec{\bff 1}}\y-\y_0\ln(\y),
\qquad
L(\x):=\y_0\Pi\x,
\qquad
\Pi:=\bbe(-A)^{-1}.
\]
Then
\[
V_{DF}(\x,\y)=R_0\,H(\y)+L(\x).
\]

We compute the two derivatives separately.

First,
\[
\dot W
=
{\vec{\bff 1}}\,\y'-\y_0\frac{\y'}{\y}
=
{\vec{\bff 1}}\bigl(\f(\y)-\diag(\y)\bbe\x\bigr)
-\y_0\frac{\f(\y)-\diag(\y)\bbe\x}{\y}.
\]
Since
\[
{\vec{\bff 1}}\diag(\y)\bbe\x=\y\,\bbe\x,
\qquad
\y_0\frac{\diag(\y)\bbe\x}{\y}=\y_0\,\bbe\x,
\]
this gives
\[
\dot W
=
{\vec{\bff 1}}\f(\y)-\y_0\frac{\f(\y)}{\y}
-\y\,\bbe\x+\y_0\,\bbe\x.
\]

Next,
\[
\dot Q
=
\y_0\Pi\x'
=
\y_0\Pi(\ba_n\y\bbe+A)\x
=
\y_0\Pi\ba_n\,\y\bbe\x+\y_0\Pi A\x.
\]
Using
\[
R_0=\y_0\Pi\ba_n,
\qquad
\Pi A=\bbe(-A)^{-1}A=-\bbe,
\]
we obtain
\[
\dot Q
=
R_0\,\y\,\bbe\x-\y_0\,\bbe\x.
\]

Therefore
\[
\dot V_{DF}
=
R_0\dot W+\dot Q
=
R_0\Bigl({\vec{\bff 1}}\f(\y)-\y_0\frac{\f(\y)}{\y}\Bigr)
+(R_0-1)\y_0\,\bbe\x.
\]
Since
\[
{\vec{\bff 1}}\f(\y)-\y_0\frac{\f(\y)}{\y}
=
(\y-\y_0)\frac{\f(\y)}{\y},
\]
and
\[
\f(\y)=-\diag(\mu_S)(\y-\y_0),
\]
it follows that
\begin{equation}\label{can}
\dot V_{DF}
=
-R_0\sum_j\frac{\mu_{S,j}(\y_j-\y_{0,j})^2}{\y_j}
+(R_0-1)\y_0\,\bbe\x.
\end{equation}

If \(R_0<1\), both terms in \eqref{can} are nonpositive, with equality only if
\[
\y=\y_0,\qquad \x=0.
\]
Hence \(V_{DF}\) is a strict Lyapunov function and the DFE is globally asymptotically stable.

If \(R_0=1\), then
\[
\dot V_{DF}
=
-\sum_j\frac{\mu_{S,j}(\y_j-\y_{0,j})^2}{\y_j}\le 0.
\]
On the set \(\{\dot V_{DF}=0\}\) one has \(\y=\y_0\), and then
\[
0=\y'=-\diag(\y_0)\bbe\x,
\]
so \(\bbe\x=0\). Therefore
\[
\x'=A\x.
\]
Since \(A\) is Hurwitz, the only invariant point is \(\x=0\). By LaSalle's invariance principle, the DFE is globally asymptotically stable also when \(R_0=1\).

If \(R_0>1\), the unstable infected block has spectral bound \(\lambda_{\max}(J_{DFE})>0\), {where $\lambda_{\max}(J_{DFE})$ is the largest real part of eigenvalues of Jacobian of infected equations with respect to infected variables, evaluated at the DFE}, so the DFE is unstable.
\end{proof}


\subsection{Sixth Law: A \Lf\ establishing the global stability of the endemic equilibrium for rank--one balanced bilinear models with one--susceptible class}\label{s:LEE}

The works \cite{IggidrCEP,Fall,Bonzi} suggested the natural question of the existence of an explicit \Lf\ which establishes the global stability of the EE of \ABa s when $R_0>1$. They offered a positive answer for the case where $A$ is bidiagonal. More recently, Earn and McCluskey \cite{EarnMc} showed that this restriction is unnecessary when $m=1$, by exploiting Volterra-type entropy identities.

\beT[GAS of the EE for \ABa s with one susceptible class]\label{thm:general_GAS}
Consider
\begin{equation}\label{eq:DSI-1S}
\begin{cases}
S'=\Lambda-\mu_S S-S\,\vbe I,\\[2pt]
I'=(S\,\ba\,\vbe + A) I,
\end{cases}
\qquad
S\in\R_+,\ I\in\R_+^n,
\end{equation}
where $\ba\in\R_+^n$ satisfies $\vec 1\,\ba=1$, $\vbe\in\R_+^{1\times n}$ is not identically zero, and $A$ is Metzler and Hurwitz.
If $R_0>1$, the endemic equilibrium $(\Se,\Ie)\in\R_{>0}\times\R_{>0}^n$ is globally asymptotically stable on $\R_+^{n+1}\setminus\{I=0\}$.
A Lyapunov function is \eqref{eq:Vshort} with weights $\vec a=\Se\,\vbe(-A)^{-1}\ge 0$.
\eeT

\begin{proof}
The endemic equilibrium is unique and satisfies
\begin{equation}\label{eq:EE-id}
\Lambda=\mu_S\Se+\Se\,\vbe\Ie,\qquad
A\Ie=-\ba\,\Se\,\vbe\Ie.
\end{equation}

Let $G(\theta)=\theta-1-\ln\theta$ and define
\begin{equation}\label{eq:Vshort}
V(S,I)=\Se\,G\!\Big(\frac{S}{\Se}\Big)+\sum_{k=1}^n a_k\,\Ie_k\,G\!\Big(\frac{I_k}{\Ie_k}\Big).
\end{equation}
Set $s=S/\Se$ and $y_k=I_k/\Ie_k$. Then
\begin{equation}\label{eq:dV-start}
\dot V=(1-\tfrac1s)S' + \sum_{k=1}^n a_k(1-\tfrac1{y_k})I_k'.
\end{equation}

\medskip\noindent
\textbf{Step 1: the $S$-term.}
Using \eqref{eq:EE-id} and the identity $(A-B)\sum_i T_i$ with zero-sum equilibrium weights, one obtains
\begin{equation}\label{eq:dVS}
(1-\tfrac1s)S'
= -\Lambda\,G(1/s)-\mu_S\Se\,G(s)
+\sum_{j=1}^n \Se\b_j\Ie_j\bigl[G(y_j)-G(y_js)\bigr].
\end{equation}

\medskip\noindent
\textbf{Step 2: the $I$-terms.}
Using \eqref{eq:EE-id},
\begin{equation}\label{eq:Ik-split}
(1-\tfrac1{y_k})I_k'
= \sum_j\ba_k \Se\b_j\Ie_j(1-\tfrac1{y_k})(y_js-y_k)
+\sum_j A_{kj}\Ie_j(1-\tfrac1{y_k})(y_j-y_k).
\end{equation}
Applying the entropy identity $(A_1-A_2)(B_1-B_2)$ to each term yields
\begin{equation}\label{eq:dVk}
(1-\tfrac1{y_k})I_k'
= \sum_j \ba_k \Se\b_j\Ie_j\bigl[G(y_js)-G(y_k)-G(y_js/y_k)\bigr]
+\sum_j A_{kj}\Ie_j\bigl[G(y_j)-G(y_k)-G(y_j/y_k)\bigr].
\end{equation}

\medskip\noindent
\textbf{Step 3: cancellation.}
Summing \eqref{eq:dVS} and $\sum_k a_k$ times \eqref{eq:dVk}, cancellation of all $G(y_\ell)$ terms requires
\begin{equation}\label{eq:AT-system}
\vec a\,A = -\Se\vbe,
\end{equation}
which has the unique solution
\begin{equation}\label{eq:a-nonneg}
\vec a=\Se\,\vbe(-A)^{-1}\ge 0.
\end{equation}
The $G(y_\ell s)$ terms cancel since $\vec a\ba=\Se\,\vbe(-A)^{-1}\ba=\Se\,\mR=1$.

\medskip

Thus
\begin{equation}\label{eq:dV-final}
\dot V
=
-\Lambda\,G(1/s)
-\mu_S\Se\,G(s)
-\sum_{i,j} c_{ij}\,G\!\Big(\frac{y_j}{y_i}\Big)
-\sum_{i,j} d_{ij}\,G\!\Big(\frac{y_j s}{y_i}\Big),
\end{equation}
with $c_{ij},d_{ij}\ge 0$. Hence $\dot V\le 0$ on $\R_{>0}^{n+1}$.

Moreover, $\dot V=0$ implies $s=1$ and $y_j=y_i$ for all $i,j$, hence $I=\Ie$.
By LaSalle's invariance principle, $(\Se,\Ie)$ is globally asymptotically stable on $\R_+^{n+1}\setminus\{I=0\}$.
\end{proof}

\beR[Relay relation between DFE and EE Lyapunov weights]\label{r:relay}
With $m=1$, the endemic equilibrium satisfies $\ye=1/\mR$
(Fourth Law, \eqref{n1}), which is C-independent.
Denoting $\vec a_0:=\y_0\bbe(-A)^{-1}$ the DFE Lyapunov weight
(Remark~\ref{r:NGMeig}), the EE Lyapunov weight $\vec a$ satisfies the
\emph{relay relation}
\[
\vec a = \ye\,\bbe(-A)^{-1}.
\]
Thus the EE and DFE Lyapunov weights are proportional, related by the
C-independent scalar $\y_0/\ye=\y_0 \mR=R_0$, and are in fact equal if their Volterra parts are assumed to have the same dependence in $S$ (but with different $\y_0,\ye$).
\eeR

\subsection{Seventh Law: A determinant law for rank--one balanced bilinear models with one--susceptible class}\lbl{s:det}

The determinant identity proved in \cite{AAH22} for the SAIR model
\[
\det(J_{EE})=-\det(J_{DFE}),
\]
extends to the whole class of rank--one balanced bilinear models with one susceptible class (in this case (P) and (B) coincide):
\begin{equation}\label{eq:1S-rank1}
\begin{cases}
\y'=\Lambda-\mu_S \y-\y\,\vbe\x,\\[2mm]
\x'=\ba\,\y\,\vbe\x + A \x,
\end{cases}
\qquad
\ba\in\mathbb R^n_+,\ \vec{\mathbf 1}\ba=1,
\end{equation}
with $A$ Metzler--Hurwitz.

\beT[Determinant law]\label{prop:detlaw}
Consider \eqref{eq:1S-rank1}, and assume $R_0>1$, where
\[
R_0=\frac{\Lambda}{\mu_S}\,\vbe(-A)^{-1}\ba.
\]
Let $J_{DFE}$ and $J_{EE}$ denote the Jacobians at the disease--free and endemic equilibria.
Then
\begin{equation}\label{eq:detlaw}
\det(J_{EE})=-\det(J_{DFE}).
\end{equation}
More explicitly,
\begin{equation}\label{eq:detexplicit}
\det(J_{DFE})=-\mu_S\det(A)\,(1-R_0),
\qquad
\det(J_{EE})=\mu_S\det(A)\,(1-R_0).
\end{equation}
\eeT

\begin{proof}
At the DFE, $\y_0=\Lambda/\mu_S$ and $\x=0$, hence
\[
J_{DFE}=
\begin{pmatrix}
-\mu_S & -\y_0\vbe\\
0 & A+\y_0\ba\vbe
\end{pmatrix},
\]
so
\[
\det(J_{DFE})=-\mu_S\,\det(A+\y_0\ba\vbe).
\]
By the matrix determinant lemma,
\[
\det(A+\y_0\ba\vbe)
=\det(A)\bigl(1+\y_0\vbe A^{-1}\ba\bigr).
\]
Since
\[
R_0=\y_0\,\vbe(-A)^{-1}\ba=-\y_0\,\vbe A^{-1}\ba,
\]
it follows that
\[
\det(A+\y_0\ba\vbe)=\det(A)(1-R_0),
\]
which yields the first identity in \eqref{eq:detexplicit}.

\medskip

At the endemic equilibrium $(\ye,\xe)$, define $q^*:=\vbe\xe$. Then
\[
J_{EE}=
\begin{pmatrix}
-\mu_S-q^* & -\ye\vbe\\
q^*\ba & A+\ye\ba\vbe
\end{pmatrix}.
\]

Perform the row operation $R_{i+1}\mapsto R_{i+1}+\ba_i R_1, i=1,...,n$ to obtain
\[
\widetilde J_{EE}=
\begin{pmatrix}
-\mu_S-q^* & -\ye\vbe\\
-\mu_S\ba & A
\end{pmatrix}.
\]
Hence, {by Schur's formula}
\[
\det(J_{EE})
=\det(A)\Bigl(-\mu_S-q^* -(-\ye\vbe)A^{-1}(-\mu_S\ba)\Bigr).
\]

Using the endemic identity
\[
A\Ie=-\ba\,\Se\,\vbe\Ie,
\]
we obtain
\[
\det(J_{EE})
=\det(A)\bigl(-\mu_S-q^*+\mu_S\bigr)
=-q^*\det(A).
\]
From $0=\Lambda-\mu_S \ye-\ye q^*$ and $\ye=\y_0/R_0$, it follows that
$q^*=\mu_S(R_0-1)$, hence
\[
\det(J_{EE})=\mu_S(1-R_0)\det(A),
\]
which proves \eqref{eq:detlaw}.
\end{proof}

\begin{remark}[Index interpretation]\label{rem:index}
Assume that all equilibria are hyperbolic.
The local index of an equilibrium is $\mathrm{sign}(\det J)$.

At the DFE,
\[
\det(J_{DFE})=-\mu_S\det(A)(1-R_0),
\]
and since $A$ is Hurwitz, $\det(A)=(-1)^n|\det(A)|$.
Thus the index of the DFE changes sign at $R_0=1$.

At the EE,
\[
\det(J_{EE})=-\det(J_{DFE}),
\]
so the two equilibria have opposite index.

This is consistent with degree theory: on any positively invariant compact region of $\R_+^{n+1}$ containing both equilibria and no others, the sum of indices is zero. Hence the appearance of the endemic equilibrium at $R_0=1$ is accompanied by a sign change of the Jacobian determinant, and the two equilibria necessarily have opposite stability types.
\end{remark}

\begin{remark}[Limitations]
The identity $\det(J_{EE})=-\det(J_{DFE})$ relies on:
\begin{enumerate}
\item the rank--one infection structure $\ba\,\y\,\vbe\x$;
\item the single susceptible equation, which allows elimination of the rank--one perturbation at the EE.
\end{enumerate}
For multi--susceptible systems, this mechanism fails and no general determinant sign law is expected.
\end{remark} 
\def\x{\boldsymbol{x}}\def\xe{\bar{\x}}%
\def\y{\boldsymbol{y}}\def\ye{\bar{\y}}%

\section{Effect of feedback $C I$ in $S'$ on the seven laws, case~(P)}\label{s:7LC}

\beT[Effect of feedback $CI$ in $S'$ on the seven laws, case~(P)]\label{thm:7laws-C}
Consider the \ABP\ \eqref{ABP} with recovery feedback $C\in\R_+^{m\times n}$, $m\ge1$.
Since $CI$ vanishes at the DFE and leaves the infection subsystem invariant,
the DFE $S_0$, $\mR=\bbe(-A)^{-1}\ba$, $R_0 =S_0 \mR$, the normalization $ \y_e\,\mR=1$, and $D_w=(-A)^{-1}\ba$ are $C$-independent and unchanged. The direction of $\bar I$~\eqref{xee}, $\bar I\propto D_w$, is also preserved.

The modifications are:

\BEN
\im \textbf{For any $m>1$ (Case~P).}
The modified Law~3:
\begin{equation}\label{eq:ye-C}
\bar S = \bigl(k\,\diag[\mR]-A_S\bigr)^{-1}(\Lambda + C\bar I),
\qquad \bar I = k D_w.
\end{equation}

The normalization condition $\vec{\bar y}\,\mR=1$ yields
\begin{equation}\label{eq:Hc}
H_C(k):=\vec{\bff 1}\,\mR\,(k\,\diag[\mR]-A_S)^{-1}(\Lambda + k C D_w)=1.
\end{equation}

If $R_0>1$, at least one positive solution of \eqref{eq:Hc} exists.

If, moreover,
\begin{equation}\label{eq:bound}
\| C D_w \|_{\infty}
<
\frac{\min_i \mu_{S,i}}{\max_i \mR_i},
\end{equation}
then \eqref{eq:Hc} admits a unique positive solution.

\im \textbf{For $m=1$,}
the system \eqref{ABP} becomes
\begin{equation}\label{eq:DSI-1S-C}
\begin{cases}
S'=\Lambda-\mu_S S-S\,\vbe I+CI,\\[2pt]
I'=\pp{S\,\ba\,\vbe+A}I,
\end{cases}
\qquad C\in\R_+^{1\times n}.
\end{equation}
Laws~0--5 and~7 are unchanged.
Only Law~6 is affected: the derivative of the EE Lyapunov function acquires the additional term
\[
\Bigl(1-\frac{\bar S}{S}\Bigr)\sum_j C_j I_j,
\]
so global asymptotic stability of the EE is no longer automatic.
\EEN
\eeT

\subsubsection*{Conjectures and Observations}

\BEN
\im \textbf{[Loss of monotonicity]:}
Let $D=\diag[\mR]$ and $M(k)=(kD-A_y)^{-1}$. Then
\[
H_C(k)=\vec{\bff 1}\,\mR\,M(k)\Lambda
+
k\,\vec{\bff 1}\,\mR\,M(k) C D_w.
\]
The first term is strictly decreasing, while the second is nonnegative but not monotone in general. Hence $H_C$ need not be monotone, in contrast with the case $C=0$.

\im \textbf{[Multiplicity of endemic equilibria]:}
When \eqref{eq:bound} fails, $H_C$ may admit several positive roots, corresponding to several endemic equilibria. Thus uniqueness of the EE is no longer structural.

\im \textbf{[Backward bifurcation mechanism]:}
Backward bifurcation may occur: \eqref{eq:Hc} may admit a positive solution even when $R_0<1$. A sufficient mechanism is the existence of $k_1,k_2>0$ such that
\[
H_C(k_1)>1,\qquad H_C(k_2)<1,
\]
which implies, by continuity, at least one solution of \eqref{eq:Hc}.

\im \textbf{[Index and determinant interpretation in the one--susceptible case]:}
In the one--susceptible case $m=1$, each positive root of \eqref{eq:Hc} yields an endemic equilibrium. Thus loss of monotonicity of $H_C$ is the scalar signature of multiplicity of endemic equilibria.

When the corresponding equilibria are hyperbolic, changes in the number of positive roots are accompanied by changes of local index, and hence by sign changes of the Jacobian determinant. In particular, saddle--node creation or annihilation of endemic equilibria corresponds to a zero of the Jacobian determinant along the endemic branch. This connects the feedback mechanism studied here with the determinant and index considerations of the previous subsection.

\im \textbf{[Scalar reduction, index, and saddle--node mechanism]:}
In the one--susceptible case $m=1$, the endemic equilibria are in one-to-one correspondence with the positive roots of the scalar equation \eqref{eq:Hc}. Thus the feedback problem reduces to the study of the real function $H_C(k)-1$.

Assume that the corresponding equilibria are hyperbolic. Then each simple root $k^*$ of \eqref{eq:Hc} defines an equilibrium whose local index is given by
\[
\operatorname{ind}(E(k^*))=\mathrm{sign}\bigl(\det J(k^*)\bigr).
\]

Loss of monotonicity of $H_C$ allows the creation or annihilation of pairs of roots. At such critical parameter values, there exists $k^*$ such that
\[
H_C(k^*)=1,\qquad H_C'(k^*)=0,
\]
which characterises a saddle--node bifurcation of endemic equilibria.

At this point, the Jacobian determinant vanishes:
\[
\det J(k^*)=0,
\]
and the index changes sign. Thus multiplicity of endemic equilibria corresponds precisely to index redistribution along the endemic branch.

This provides a direct connection between:
\begin{itemize}
\item the scalar equation \eqref{eq:Hc},
\item the loss of monotonicity induced by the feedback $C$,
\item and the determinant sign structure described in Section~\ref{s:det}.
\end{itemize}

In particular, backward bifurcation corresponds to the appearance of a pair of equilibria of opposite index below threshold.

\EEN

\section{The computation of Perron eigenvectors
for matrices admitting regular splittings in the sense of \cite{Varga}, via Kirchhoff's matrix-tree theorem}
\label{s:KP}

In this section, we recall and apply the computation of Perron eigenvectors
 via Kirchhoff's matrix-tree theorem.

We provide first a Kirchhoff matrix-tree theorem \ref{thm:KP}  which applies both to computing stationary distributions of Markov processes and to computing Lyapunov weights at the DFE for ME (or CRN) models, under a condition of \cite{Shuai13}.

We consider Metzler matrices $J$ admitting a regular splitting \cite{Varga}
\[
    J = F - V,
\]
where $F\ge0$ and $V$ is a nonsingular $M$-matrix with $V^{-1}\ge0$ \cite{Varga}.
When $J$ arises as the transversal Jacobian of a compartmental epidemic model at the DFE (with tangential variables $\y$ fixed), it describes the linearised infection dynamics.

If the matrix $J$  is also diagonally dominant (with nonnegative off-diagonal entries and nonpositive row sums), then it defines,  after transposition if needed, the generator of a continuous-time Markov process with possible killing.
Below threshold ($\lambda_P<0$), it is a subgenerator, and  at threshold ($\lambda_P=0$) it becomes conservative.
This interpretation connects Markov processes with killing (see for example \cite{Norris97}) to epidemic models \cite{Kendall,Diek,Van,Shuai13}.

 {\bf For a conservative irreducible generator $Q$}, Kirchhoff's Markov matrix-tree theorem states that the unique (up to scaling) positive left null vector $\vec \pi>0$ satisfying $\vec \pi Q=0$ has components proportional to the {\bf diagonal cofactors} $C_{kk}(-Q)$ \cite{Norris97,GuoLiShuai06}. Here
\[
C_{kk}(M)=\det M^{(k,k)},
\]
where $M^{(k,k)}$ is obtained from $M$ by deleting the $k$-th row and the $k$-th column.
Also, for a conservative irreducible generator $Q$, rank$(Q)=n-1$, and
$\operatorname{adj}(Q)=\{(-1)^{k+j}C_{kj}(Q)\}_{j,k}$ is of rank one, with rows proportional to $\vec\pi$ and columns proportional to $\mathbf{1}$.

In our setting, the role of  (sub)generator is played by $J=F-V$,   while $V-F$ is the associated $M$-matrix entering the adjugate formula.

\beT[Kirchhoff--Perron formula for regular splittings]
\label{thm:KP}
Let $J$ be an irreducible Metzler matrix with spectral abscissa $\lambda_P$. Then $\lambda_P$ is a simple eigenvalue of $J$.

\begin{enumerate}
\item \emph{(Perron--Frobenius.)}
There exist unique (up to scaling) positive left and right Perron eigenvectors:
\[
   \vec \pi \, J = \lambda_P\, \vec \pi ,\qquad J\,\bw  = \lambda_P\,\bw ,
   \qquad \vec \pi >0,\quad \bw >0.
\]

\item \emph{(Adjugate factorisation at the Perron root.)}
The matrix $J-\lambda_P I$ is singular with rank $n-1$, and its adjugate has rank one and factors as
\[
   \operatorname{adj}(J-\lambda_P I)= c\,\bw \,\vec \pi ,\qquad c\neq 0.
\]
Hence all columns of $\operatorname{adj}(J-\lambda_P I)$ are proportional to $\bw$, all rows are proportional to $\vec \pi$, and
\[
   C_{kk}(\lambda_P I-J)\propto \bw_k\,\pi_k.
\]
\end{enumerate}
\end{theorem}

\begin{remark}
In the Markov case ($\lambda_P=0$), $J$ is a conservative generator ($J\mathbf{1}=0$)
and Kirchhoff's theorem applied to the  Laplacian $-J$ yields
\[
C_{kk}(-J)=c_0\,\pi_k,\qquad c_0>0.
\]
\end{remark}

\beT[Roles of NGM eigenvectors at a unique DFE]\label{cor:NGM}
Assume that the disease-free system
\[
S'=g(S,0)
\]
has a unique equilibrium $S_0>0$, where the transversal Jacobian $J=J_x(S_0,0)$ is 
irreducible Metzler and admits a  \regS:
\[
J = F(S_0) - V(S_0).
\]
 
 Define
\[
K =K(S_0):= F(S_0)V(S_0)^{-1},\qquad \T K:=V(S_0)^{-1}F(S_0),\qquad R_0 := \rho(K)=\rho(\T K).
\]
Then:
\begin{enumerate}
\item $\lambda_P(J)=0$ iff $R_0=1$.

\item Let $\bw$ denote the right Perron eigenvector of $K=F(S_0)V^{-1}$. Then
\[
D_w:=V^{-1}\bw
\]
is the right Perron eigenvector of $\T K=V^{-1}F(S_0)$, 
and
\[
J(S_0)D_w=(R_0(S_0)-1)\bw
\]
holds. 

At threshold $R_0(S_0)=1$, the vector $D_w$ is a right null vector of $J(S_0)$, and 
if an endemic branch bifurcates from $(S_0,0)$ at threshold, then its infected component is tangent to the ray spanned by $D_w$.

In rank-one Case~(P), one has moreover
\[
\bar I=k D_w,\qquad k>0,
\]
for the unique endemic equilibrium, so the whole local endemic branch lies on that ray.

\item Let $\vec \pi$ be the left Perron eigenvector of $K=F(S_0)V^{-1}$,
\[
\vec\pi K=\vec\pi F(S_0)V^{-1}=R_0\vec\pi,
\]
and let
\[
\vec q=\vec\pi V
\]
denote the left Perron eigenvector of $\T K=V^{-1}F(S_0)$. Put
\[
Q(I)=\vec q\,I=\vec\pi V I.
\]
If, moreover,
\[
I'=(F(S_0)-V)I-f(S,I),
\qquad f(S,I)\ge0,
\]
then
\[
Q'=(R_0-1)\vec\pi\,I-\vec\pi Vf(S,I),
\]
and therefore $Q$ is a transversal Lyapunov function for the invariant manifold $\{I=0\}$ whenever $R_0<1$.
\end{enumerate}
\eeT

\begin{proof}
Item~(1) is the standard equivalence for regular splittings:
\[
\lambda_P<0 \iff \rho(F(S_0)V^{-1})<1,\qquad
\lambda_P=0 \iff \rho(F(S_0)V^{-1})=1.
\]

For item~(2), let $\bw>0$ satisfy
\[
K\bw=F(S_0)V^{-1}\bw=R_0\bw,
\]
and define
\[
D_w:=V^{-1}\bw.
\]
Since $V$ is nonsingular, $D_w\neq0$. Multiplying the eigenvalue relation by $V^{-1}$ on the left gives
\[
V^{-1}F(S_0)V^{-1}\bw=R_0V^{-1}\bw,
\Eq
\T K D_w=V^{-1}F(S_0)\,D_w=R_0D_w.
\]
Hence $D_w$ is the right Perron eigenvector of $\T K$.

At threshold, $R_0=1$, and therefore
\[
J D_w=(F(S_0)-V)D_w
      =F(S_0)V^{-1}\bw-\bw
      =(K-I)\bw
      =0.
\]
Thus $D_w$ is the right null vector of $J$.

Assume now that an endemic branch
\[
(\bar S(\eps),\bar I(\eps)),\qquad \eps\ge0,
\]
bifurcates from the DFE $(S_0,0)$ at $R_0=1$, is differentiable at $\eps=0$, and satisfies
\[
(\bar S(0),\bar I(0))=(S_0,0).
\]
Write
\[
\bar I(\eps)=\eps d+o(\eps),\qquad \eps\downarrow0.
\]
Since $(\bar S(\eps),\bar I(\eps))$ is an equilibrium, the infected equilibrium equation has the form
\[
0=(F(\bar S(\eps))-V)\bar I(\eps)+N(\bar S(\eps),\bar I(\eps)),
\]
where $N(S,I)=o(\|I\|)$ as $I\to0$. Dividing by $\eps$ and letting $\eps\downarrow0$ yields
\[
(F(S_0)-V)d=Jd=0.
\]
Hence $d\in\ker J$. Since $J$ is irreducible Metzler and $0$ is its simple Perron eigenvalue at threshold, $\ker J$ is one-dimensional. Therefore
\[
d=\eta D_w
\]
for some scalar $\eta\neq0$. Thus the infected component of the endemic branch is tangent at $(S_0,0)$ to the ray spanned by $D_w$.

In rank-one Case~(P), one has
\[
F(S)=\ba\,\vec w(y),
\]
so at the unique endemic equilibrium on the branch,
\[
0=(F(\bar S)+A)\bar I=\ba(\vec w(y_e)\bar I)+A\bar I.
\]
Writing
\[
k:=\vec w(\bar S)\bar I,
\]
we obtain
\[
A\bar I=-k\ba,
\qquad\text{hence}\qquad
\bar I=k(-A)^{-1}\ba=kD_w.
\]
Thus, in rank-one Case~(P), the whole local endemic branch lies on the ray spanned by $D_w$.

For item~(3), since
\[
\vec\pi K=\vec\pi F(S_0)V^{-1}=R_0\vec\pi,
\]
right-multiplying first by $V$ yields
\[
\vec\pi F(S_0)=R_0\vec\pi V,
\]
and then right-multiplying by $V^{-1}$ gives
\[
\vec\pi V^{-1}F(S_0)=R_0\vec\pi.
\]

Thus, define
\[
\vec{q} := \vec{\pi} V,
\]
which is the left Perron eigenvector of $\tilde K = V^{-1}F(S_0)$.

Define
\[
Q(I) := \vec{q}\, I.
\]

Using the decomposition 
$
I' = (F(S_0) - V)I - f(S,I),
$, 
we compute 
$ 
Q' = \vec{q}\, I' = \vec{q}(F(S_0)-V)I - \vec{q}\, f(S,I).$

Since $\vec{q} V^{-1}F(S_0) = R_0 \vec{q}$, it follows that 
$
\vec{q}(F(S_0)-V) = (R_0 - 1)\vec{q}.
$

Therefore,
$
Q' = (R_0 - 1)\,\vec{q} I - \vec{q}\, f(S,I).
$

Because $\vec{q} > 0$, $I \ge 0$, and $f(S,I) \ge 0$, we obtain
\[
Q' \le 0 \quad \text{whenever } R_0 < 1.
\]

Therefore $Q$ is a transversal Lyapunov function for the invariant manifold $\{I=0\}$.
\end{proof}

\section{Conclusion}\lbl{s:con}

We reviewed a class of mathematical epidemiology (ME) models that we call \ABP s, for which explicit formulas for the interior fixed point and the basic reproduction number, together with elegant Lyapunov functions establishing the asymptotic stability of the disease-free equilibrium (DFE) and endemic equilibrium (EE), are available when certain matrices have rank one.

Several themes emerge. 

First, rank one plays two rather different roles. On the one hand, it yields explicit formulas for Perron eigenvectors, endemic equilibria, and Goh--Volterra coefficients. On the other hand, the DFE theory itself is not inherently rank-one: the Perron Lyapunov weight persists for higher-rank NGMs, although it is no longer explicit. Thus rank one is best viewed not as the source of the theory, but as the regime in which the general theory becomes closed-form.

Second, the relation between DFE and EE stability is stronger than is usually emphasized. In the rank-one \ABP\ setting, both rely on the same building blocks: the Perron data of the NGM, the dwell-time vector, and a normalization law. This suggests that DFE and EE Lyapunov theories should not be developed separately, but rather as two parts of a common ``stability relay''. At present this relay is fully understood only in special cases.

Third, the comparison with Shuai--van den Driessche shows a sharp distinction between what remains valid at general rank and what does not. The DFE Lyapunov weight remains Perron-theoretic, but EE Lyapunov weights are generally no longer Perron eigenvectors; they are instead related to Kirchhoff-type cofactors. Thus one should expect a genuine bifurcation in the theory at higher rank: Perron theory controls the DFE, while graph-theoretic constructions control the EE.


These observations lead to the following open problems.

\beO [General Lyapunov functions for bilinear ME models]
Is it possible to extend both \Lf s in Sections~\ref{s:LDFE} and~\ref{s:LEE} to a wider class that includes both
\ABP\ and multi-strain models, including models whose NGM has rank greater than one?
In particular, can the Lyapunov functions of \cite{Rahman,Bulh} be obtained from some unifying general formulas?

A first subproblem is to determine whether the DFE Perron term and the EE Goh--Volterra term can be completed into a single general construction, valid beyond rank one. A second is to decide whether the higher-rank EE coefficients should be sought through Perron theory, through Kirchhoff cofactors, or through a genuinely new mechanism.
\eeO

\beO[General theory for multistrain models]
Do  formulas \cite{AAN} like
\be{mys} \bc R_i =R_i(E_0), R_0=\max\pp{ R_1, R_2},\\\mR_{\T i}^j =R_j(E_i),
\\R_2 <R_{2,c}:=\fr{R_2(E_0)}{R_2(E_1)} \Lra E_1 \text{ is LAS }\ec
\ee
where $R_0,R_i$ denote the reproduction numbers, $\mR_{\T i}^j$ the invasion numbers of resident strain $i$ by strain $j$, and $R_i(x)$ the reproduction functions,
hold for some general class of epidemic models?

More generally, one may ask whether reproduction functions admit a structural interpretation analogous to that of NGMs, and whether invasion criteria can be derived from Jacobian factorizations without requiring explicit knowledge of the competing endemic equilibria.
\eeO

\beO [Structural characterization of balanced bilinear epidemic models within CRN theory]
Characterize \ABP s within CRN theory, and investigate extensions of the bilinear results under rates which are admissible in the sense of \cite{AdLS}.

This includes at least three natural questions: how to recognize \ABP\ structure directly from a reaction network; which parts of the seven-law theory survive under admissible non-bilinear rates; and whether siphons, minimal siphons, or related CRN objects provide the correct structural substitute for the rank-one assumptions.
\eeO

\beO[Persistence and boundary structure beyond the single-siphon case]
Extend the persistence part of Shuai--van den Driessche's theory from the case of a single minimal siphon to models with several minimal siphons.

A natural conjecture is that persistence above threshold should hold whenever every minimal siphon supports a repelling boundary face, equivalently whenever the restriction of the infection dynamics to each such face has spectral radius greater than one. Formulating and proving such a criterion would provide a boundary counterpart to the higher-rank DFE theory.
\eeO

\beO[Determinant laws, relay identities, and inheritance beyond rank one]
Determine whether the determinant law, relay relations between DFE and EE weights, and related inheritance phenomena extend beyond the one-susceptible rank-one setting.

Even partial extensions would clarify which features are genuinely tied to rank one and which are manifestations of a broader algebraic structure underlying epidemic Jacobians.
\eeO


\section{ Supplementary Notation and Remarks }

\medskip
\noindent\textbf{Siphon lattice and invariant faces.}
\begin{itemize}
\item $\mathcal L$ denotes the lattice generated by the minimal siphons.
\item For $\Sigma\in\mathcal L$,
\[
\mathcal F_\Sigma=\{z_i=0:\ i\in\Sigma\}
\]
is the associated forward-invariant face.
\item $E_\Sigma\in\mathrm{relint}(\mathcal F_\Sigma)$ denotes the equilibrium on $\mathcal F_\Sigma$ under consideration.
\item At $E_\Sigma$, one has
\[
z_i^*=0\ \text{for } i\in\Sigma,
\qquad
z_i^*>0\ \text{for } i\notin\Sigma.
\]
\end{itemize}

\medskip
\noindent\textbf{Face-dependent splitting and Jacobians.}
For each $\Sigma\in\mathcal L$, variables are split as
\[
z=(x_\Sigma,y_\Sigma),
\]
where $x_\Sigma=(z_i)_{i\in\Sigma}$ are the transversal (invader) variables and
$y_\Sigma=(z_i)_{i\notin\Sigma}$ are the tangential (resident) variables.

At $E_\Sigma$, the Jacobian has block triangular form
\[
Df(E_\Sigma)=
\begin{pmatrix}
J_\Sigma^\perp & 0\\
* & J_{\mathrm{tan}}(\Sigma)
\end{pmatrix},
\]
after ordering variables as $(x_\Sigma,y_\Sigma)$, where
\[
J_\Sigma^\perp:=D_{x_\Sigma}f_{x_\Sigma}(E_\Sigma),
\qquad
J_{\mathrm{tan}}(\Sigma):=D_{y_\Sigma}f_{y_\Sigma}(E_\Sigma).
\]
Thus:
\begin{itemize}
\item $J_\Sigma^\perp$ governs invasion away from $\mathcal F_\Sigma$;
\item $J_{\mathrm{tan}}(\Sigma)$ governs the dynamics on $\mathcal F_\Sigma$.
\end{itemize}

\beR
The notation $J_\Sigma^\perp$ is reserved for the full transversal Jacobian on the face $\mathcal F_\Sigma$.
If $\sigma\subseteq\Sigma$ is a released invader block, then
$J_\sigma(E_\Sigma)$ denotes the corresponding principal transversal block.
\eeR
\beP[Block triangular Jacobian at an invariant face]
Let
\[
z=(x_\Sigma,y_\Sigma)\in \R_{\ge 0}^{p}\times \R_{\ge 0}^{q},
\qquad
z' = f(z)=\bigl(f_{x_\Sigma}(x_\Sigma,y_\Sigma),\,f_{y_\Sigma}(x_\Sigma,y_\Sigma)\bigr),
\]
with $f\in C^1$ on a neighbourhood of the face
\[
\mathcal F_\Sigma:=\{x_\Sigma=0\}.
\]
Assume that $\mathcal F_\Sigma$ is forward invariant, equivalently
\[
f_{x_\Sigma}(0,y_\Sigma)\equiv 0
\qquad\text{for all }y_\Sigma\ge 0.
\]
Then, at every equilibrium
\[
E_\Sigma=(0,y_\Sigma^*)\in \mathcal F_\Sigma,
\]
the Jacobian has block triangular form
\[
Df(E_\Sigma)=
\begin{pmatrix}
J_\Sigma^\perp & 0\\
* & J_{\mathrm{tan}}(\Sigma)
\end{pmatrix},
\]
where
\[
J_\Sigma^\perp=D_{x_\Sigma}f_{x_\Sigma}(E_\Sigma),
\qquad
J_{\mathrm{tan}}(\Sigma)=D_{y_\Sigma}f_{y_\Sigma}(E_\Sigma).
\]
\eeP

\begin{proof}
Since $f$ is $C^1$, its Jacobian at $E_\Sigma$ may be written in block form
\[
Df(E_\Sigma)=
\begin{pmatrix}
D_{x_\Sigma}f_{x_\Sigma}(E_\Sigma) & D_{y_\Sigma}f_{x_\Sigma}(E_\Sigma)\\
D_{x_\Sigma}f_{y_\Sigma}(E_\Sigma) & D_{y_\Sigma}f_{y_\Sigma}(E_\Sigma)
\end{pmatrix}.
\]
Thus it remains only to prove that
\[
D_{y_\Sigma}f_{x_\Sigma}(E_\Sigma)=0.
\]

By invariance of the face, one has
\[
f_{x_\Sigma}(0,y_\Sigma)\equiv 0
\qquad\text{for all }y_\Sigma\ge 0.
\]
Therefore the map
\[
y_\Sigma \longmapsto f_{x_\Sigma}(0,y_\Sigma)
\]
is identically zero on a neighbourhood of $y_\Sigma^*$. Differentiating with
respect to $y_\Sigma$ at $y_\Sigma=y_\Sigma^*$ gives
\[
D_{y_\Sigma}f_{x_\Sigma}(0,y_\Sigma^*)=0.
\]
Hence
\[
Df(E_\Sigma)=
\begin{pmatrix}
D_{x_\Sigma}f_{x_\Sigma}(E_\Sigma) & 0\\
D_{x_\Sigma}f_{y_\Sigma}(E_\Sigma) & D_{y_\Sigma}f_{y_\Sigma}(E_\Sigma)
\end{pmatrix},
\]
which is exactly
\[
Df(E_\Sigma)=
\begin{pmatrix}
J_\Sigma^\perp & 0\\
* & J_{\mathrm{tan}}(\Sigma)
\end{pmatrix}.
\]

\end{proof}

\medskip
\noindent\textbf{Regular splitting and face next-generation matrices.}
For an invader block $\sigma\subseteq\Sigma$,
\begin{itemize}
\item
\[
J_\sigma(E_\Sigma)=F_\sigma(E_\Sigma)-V_\sigma(E_\Sigma)
\]
denotes a regular splitting, with $F_\sigma\ge 0$ and $V_\sigma^{-1}\ge 0$;
\item
\[
K_\sigma(E_\Sigma)=F_\sigma(E_\Sigma)V_\sigma(E_\Sigma)^{-1}
\]
is the face next-generation matrix;
\item
\[
\widetilde K_\sigma(E_\Sigma)=V_\sigma(E_\Sigma)^{-1}F_\sigma(E_\Sigma)
\]
is the cyclically permuted next-generation matrix;
\item
\[
R_\sigma(E_\Sigma)=\rho\bigl(K_\sigma(E_\Sigma)\bigr)
\]
is the face reproduction number.
\end{itemize}
The classical basic reproduction number $R_0$ is the special case corresponding to the DFE face.

\medskip
\noindent\textbf{Spectral abscissae and Perron invasion directions.}
\begin{itemize}
\item
\[
\alpha_\sigma(E_\Sigma):=s\bigl(J_\sigma(E_\Sigma)\bigr)
\]
is the spectral abscissa of the released invader block $\sigma$ at $E_\Sigma$.
\item $w_\sigma\gg0$ denotes a right Perron eigenvector of $K_\sigma(E_\Sigma)$.
\item $\tws\gg0$ denotes a right Perron eigenvector of $\widetilde K_\sigma(E_\Sigma)$.
\item
\[
D_{w,\sigma}:=V_\sigma(E_\Sigma)^{-1}w_\sigma
\]
is the Perron invasion direction.
\item $a_\sigma>0$ denotes the amplitude scalar fixing the equilibrium on the released face.
\end{itemize}

\medskip
\noindent\textbf{Face-adapted Lyapunov functions.}
For $\Sigma\in\mathcal L$,
\[
V_\Sigma(z)=
\sum_{i\notin\Sigma}a_i^\Sigma\,z_i^*\,
G\!\Bigl(\frac{z_i}{z_i^*}\Bigr)
+
\sum_{j\in\Sigma}b_j^\Sigma\,z_j,
\]
with coefficients
\[
a_i^\Sigma>0,\qquad b_j^\Sigma>0.
\]
Here:
\begin{itemize}
\item $a_i^\Sigma$ are the Volterra (resident) weights;
\item $b_j^\Sigma$ are the linear (invader) weights.
\end{itemize}
For rank-one \ABa\ with $m=1$, one has the relay relation
\[
\vec l=\bar S\,\vec\pi.
\]

\medskip
\noindent\textbf{Flux, Laplacian, and Kirchhoff notation.}
\begin{itemize}
\item $\Phi_{ij}(z)\ge 0$ denotes the directed flux from species $i$ to species $j$ in a decomposition
\[
f_i=\sum_k(\Phi_{ki}-\Phi_{ik}).
\]
\item
\[
c_{ij}:=\partial_{z_i}\Phi_{ij}(z^*)\ge 0
\]
denotes the linearized flux coefficient at equilibrium.
\item $\kappa_{\Sigma\Sigma'}\ge 0$ denotes the edge weight in the relay graph associated with the inclusion or release step under consideration.
\item $L$ denotes the Laplacian of the relay graph:
\[
L_{\Sigma\Sigma}=\sum_{\Sigma'\sim\Sigma}\kappa_{\Sigma\Sigma'},
\qquad
L_{\Sigma\Sigma'}=-\kappa_{\Sigma\Sigma'}\ \text{for }\Sigma'\neq\Sigma.
\]
\item
\[
La=B
\]
denotes the Laplacian balance system for the resident weights.
\item $C_{\Sigma\Sigma}(L)$ denotes the diagonal cofactor of $L$.
By Kirchhoff's matrix-tree theorem, one has
\[
a_\Sigma\propto C_{\Sigma\Sigma}(L)
\]
when $La=0$.
\end{itemize}

\medskip
\noindent\textbf{CRN notation.}
The stoichiometric matrix is $\Gamma$. A reaction  $\rho$ is the index of a column in
$\Gamma$, and \[
\gamma_\rho:=O(\rho)-S(\rho)\in\Z^n
\]
is the stoichiometric vector of $\rho$, decomposed as a difference between output and and  source  \nne\ vectors. In the mass action case, $\rho$ may be viewed as a transition between two complexes $y\to y'$.

A reaction may be viewed also as a pair
\[
\rho=\bigl(\Gam(\rho),\Gap(\rho)\bigr),
\]
where $\Gam(\rho)$ is the set of incoming species with negative signs and $\Gap(\rho)$ is the set of outcoming species with positive signs.

 $r_\rho(z)\ge 0$ is the rate of $\rho$;
under mass action,
\[
r_\rho(z)=\kappa_\rho\,z^{S(\rho)}.
\]

\beR
The symbol $\rho$ has two distinct meanings:
as a reaction index/ hyper-arc  in expressions such as
$r_\rho$, $\gamma_\rho$, $\sum_\rho$,
and as the spectral radius in expressions such as $\rho(K)$.
The distinction is always by context.
\eeR

\begin{center}
  \scalebox{0.85}{
  \begin{tabular}{|c|c|c|}
  \hline
  Object & Role in dynamics & Transformations \\
  \hline
  $\tilde w$ (right evec of $\T K$)
  & Escape direction
  & $w:=(-A)^{-1}\tilde w$  \\
  \hline
  $\tilde\pi^T$ (left evec of $\T K$)
  & $-$
  & $\pi^T:=\tilde\pi^T(-A)^{-1}$  \\
  \hline
  $w$ (right evec of $K$)
  & $-$
  & $\tilde w=(-A)\,w$ \\
  \hline
  $\pi^T$ (left evec of $K$)
  & Lyapunov weights
  & $\tilde\pi^T=\pi^T(-A)$ \\
  \hline
  \end{tabular}
  }
  \end{center}


\bibliographystyle{plain}
\bibliography{references}

\end{document}